\documentclass[12pt]{amsart}

\input xy
\xyoption{all}

\usepackage{geometry}
\usepackage{amsmath}
\usepackage{amssymb}
\usepackage{cite}
\usepackage{amsthm}  % see geometry.pdf on how to lay out the page. There's lots.
\newtheorem{same}{This should never appear}[section]
\newtheorem{defin}[same]{Definition}
\newtheorem{claim}[same]{Claim}
\newtheorem{remark}[same]{Remark}
\newtheorem{theorem}[same]{Theorem}

\newtheorem{lemma}[same]{Lemma}
\newtheorem{fact}[same]{Fact}

\newtheorem{cor}[same]{Corollary}
\newtheorem{prop}[same]{Proposition}

\newtheorem{hypothesis}[same]{Hypothesis}

\newtheorem{defin*}{Definition}
\newtheorem*{theorem*}{Theorem}
\newbox\noforkbox \newdimen\forklinewidth
\setbox0\hbox{$\textstyle\smile$}
\setbox1\hbox to \wd0{\hfil\vrule width \forklinewidth depth-2pt
 height 10pt \hfil}
\setbox\noforkbox\hbox{\lower 2pt\box1\lower 2pt\box0\relax}
\def\unionstick{\mathop{\copy\noforkbox}\limits}

\def\nonfork_#1{\unionstick_{\textstyle #1}}

\setbox0\hbox{$\textstyle\smile$}
\setbox1\hbox to \wd0{\hfil{\sl /\/}\hfil}
\setbox2\hbox to \wd0{\hfil\vrule height 10pt depth -2pt width
              \forklinewidth\hfil}
\newbox\doesforkbox
\setbox\doesforkbox\hbox{\lower 2pt\box1 \lower 2pt\box2\lower2pt\box0\relax}
\def\nunionstick{\mathop{\copy\doesforkbox}\limits}

\def\fork_#1{\nunionstick_{\textstyle #1}}

\newcommand{\sea}{\mathfrak{C}}

\newcommand{\Mod}{\textrm{Mod}}

\newcommand{\cf}{\text{cf }}
\newcommand{\rest}{\upharpoonright}

\newcommand{\dnf}{\unionstick}
\newcommand{\dnks}{{}^{S}\dnf}

\newcommand{\bigM}{\widehat{M}}
\newcommand{\bigN}{\widehat{N}}

\newcommand{\adnf}{\dnf^*}

\newcommand{\cont}{\textrm{Contents}}

\newcommand{\seq}[1]{\langle #1 \rangle}
\newcommand{\te}[1]{\text{#1}}

\title{Forking in Short and Tame Abstract Elementary Classes}
\date{\today\\
Part of this material is based upon work done while the first author was supported by the National Science Foundation under Grant No. DMS-1402191. 
} % delete this line to display the current date

\author{Will Boney}
\email{wboney@math.harvard.edu}
\address{Mathematics Department \\ Harvard University\\ Cambridge, Massachusetts, USA}

\author{Rami Grossberg}
\email{rami@cmu.edu}
\address{Department of Mathematical Sciences \\ Carnegie Mellon University \\ Pittsburgh, Pennsylvania, USA}

\begin{document}

\maketitle

%%%%%%%%%%%%%%%%%%%%%%%%%%%%%%%%%%%%%%%%%%%%%%%
\begin{abstract}
We develop a notion of forking for Galois-types in the context of Abstract Elementary Classes (AECs).  Under the hypotheses that an AEC $K$ is tame, type-short, and failure of an order-property, we consider

\begin{defin*}
Let $M_0 \prec N$ be models from $K$ and $A$ be a set.  We say that the Galois-type of $A$ over $N$ \emph{does not fork over $M_0$}, written $A \dnf_{M_0} N$, iff for all small $a \in A$ and all small $N^- \prec N$, we have that Galois-type of $a$ over $N^-$ is realized in $M_0$.
\end{defin*}

Assuming property $(E)$ (Existence and Extension, see Definition 3.3) we show that this non-forking is a well behaved notion of independence, in particular satisfies symmetry and uniqueness and has a corresponding U-rank.  We find conditions for a universal local character, in particular derive superstability-like property from little more than categoricity in a ``big cardinal''.  Finally, we show that under large cardinal axioms the proofs are simpler and the non-forking is more powerful.

In \cite{bgkv}, it is established that, if this notion is an independence notion, then it is the only one.
\end{abstract}
%%%%%%%%%%%%%%%%%%%%%%%%%%%%%%%%%%%%%%%%%%%%%%%

\tableofcontents

%%%%%%%%%%%%%%%%%%%%%%%%%%%%%%%%%%%%%%%%%%%%%%%
\section{Introduction}
%%%%%%%%%%%%%%%%%%%%%%%%%%%%%%%%%%%%%%%%%%%%%%%

Much of modern model theory has focused on Shelah's forking.  In the last twenty years, significant progress has been made towards understanding of unstable theories, especially simple theories (Kim \cite{Kim} and Kim and Pillay \cite{KP1}), $NIP$ theories (surveys by Adler \cite{Adler} and Simon \cite{Simon}), and, most recently, $NTP_2$ (Ben-Yaacov and Chernikov \cite{Ch1} and Chernikov, Kaplan, and Shelah \cite{Ch2}). 

  In the work on classification theory for Abstract Elementary Classes (AECs), such a nicely behaved notion is not known to exist.  However, much work has been done towards this goal.  %Already in 1970, Shelah \cite{Sh3} introduced splitting as a weak independence notion for a nonelementary context that is now known as homogeneous model theory (see Grossberg and Lessmann \cite{GrLe1} and Buechler and Lessmann \cite{BuLe}).  For the more general cases of classes axiomatizeable by an $L_{\lambda^+,\omega}$ sentence or AECs, very little is known in this direction, although there have been several attempts.
Around 2005, homogeneous model theory--working under the assumption that there exists a monster model which is sequential homogeneous (but not necessarily saturated as in the first-order sense) and types consists of sets of first-order formulas--reached a stage of development that parallels that of first-order model theory in the seventies.  There is a Morley-like categoricity theorem (Keisler \cite{KeLomega} and Lessmann \cite{lessman2000}), forking exists (Buechler and Lessman \cite{BuLe}), and even a main gap is true (Grossberg and Lessmann  \cite{GrLe05}).  Hyttinen and Kes{\"a}l{\"a} studied a further extension of homogeneous model theory called finitary AECs in \cite{hyttinenkesala2006} and in \cite{hyttinenkesala11}.  They established both Morley’s categoricity theorem and that non-splitting is a variant of forking under the assumptions of $\aleph_0$-stability and what they call simplicity (like our extension property) in a countable language.

However, as AECs are much more general the situation for AECs is  more  complicated.  There are classes axiomatized by $L_{\omega_1, \omega}$ that do not fit into the framework of homogeneous model theory:
\begin{enumerate}
	\item Marcus \cite{marcus75} constructed an $L_{\omega_1,\omega}$ sentence that is categorical in all cardinals but does not have even an $\aleph_1$-homogenous model.
	\item Hart and Shelah \cite{hash323} constructed, for each $k < \omega$, an $L_{\omega_1,\omega}$ sentence $\psi_k$ which is categorical in all $\aleph_n$ for $n\leq k$ but not categorical in higher cardinals.  By the categoricity theorem for finitary AECs \cite{hyttinenkesala11},  this means that $\Mod(\psi_k)$ is not homogeneous as it is not even finitary.
\end{enumerate}

In \cite[Chapter N]{shelahaecbook}, Shelah explains the importance of classification theory for AECs.  At the referee's suggestion, we summarize the argument here, although the truly interested reader should consult the source.

As mentioned above, classification theory has become the main focus of model theory.  Shelah and other early workers were motivated by purely abstract problems, such as the main gap in \cite{shelahfobook}.  The machinery used to solve these problems turned out to be very powerful and, about 20 years later, Chatzidikis, Hrushovski, Scanlon, and others discovered deep applications to geometry, algebra, and other fields.

However, this powerful machinery was restricted because it only applied to first-order model theory.  This is natural from a logical point of view as first-order logic has many unique features (e. g., compactness), but there are many mathematical classes that are not first order axiomatizable: we list some in this introduction and in Section \ref{fullsitsection} and each of Grossberg \cite{grossberg2002}, Baldwin \cite{baldwinbook}, and \cite[Chapter N]{shelahaecbook} contain their own lists.  The logic needed to axiomatize each context varies, from $L_{\omega_1, \omega}(Q)$ for quasiminimal classes to $L_{|R|^+, \omega}$ for torsion $R$-modules. Varying the substructure relation (e. g., subgroup vs. pure subgroup) complicates the picture further.

A unifying perspective is given by AECs and Shelah began their classification (and their study) in the late 1970s.  Here again questions of number of nonisomorphic models have formed the basic test questions.  The most central one here is Shelah's Categoricity Conjecture; Shelah proposed this conjecture  for $L_{\omega_1, \omega}$ in the late seventies as a way to measure the development of the relevant classification theory.  At present, there are many partial results that approximates this conjecture and harder questions for AECs.  Despite an estimated of more than 2,000 published pages, the full conjecture is not within reach of current methods, in contrast to the existence of relatively simple proofs of the conjecture for the cases of homogeneous models and finitary AECs.  Due to the lack of compactness and syntax, extra set-theoretic assumptions (in addition to new techniques) have been needed to get these results; the strong Devlin-Shelah diamonds on successors in \cite{sh576} and large cardinals in \cite{makkaishelah} are excellent examples of this.

Differing from Shelah, our vision is that model-theoretic assumptions (especially tameness and type-shortness here) will take the place of set-theoretic ones.  The hope here is two-fold: first, that, although these assumptions don't hold everywhere, they can be shown to hold in many natural and, second, that these model-theoretic assumptions are enough to develop a robust classification theory.  This paper (and follow-ups by Boney, Grossberg, Kolesnikov, and Vasey) provide evidence for the first hope and examples described in Section \ref{fullsitsection} provide evidence for the first.

In \cite{sh394}, Shelah introduced analogues of splitting and strong splitting for AECs.  Building on this, Shelah \cite{sh576}, and Grossberg, VanDieren, and Villaveces \cite{gvv}, Grossberg and VanDieren \cite{tamenesstwo, tamenessthree} used tameness to prove an upward categoricity transfer.  This helped cement tameness as an important property in the classification of AECs.  Working in a stronger context, Makkai and Shelah \cite{makkaishelah} studied the case when a class is axiomatized by an $L_{\kappa,\omega}$ theory and $\kappa$ is strongly compact.   They managed to obtain an eventual categoricity theorem by introducing a forking-like relation on types.  In this particular case, Galois types (defined in \S \ref{prelimsection}) can be identified with complete set of formulas taken from a fragment of $L_{\kappa,\kappa}$.  In their paper, Makkai and Shelah assumed not only that $\kappa$ is strongly compact but also that the class of structures is categorical in some $\lambda^+$ where $\lambda\geq\beth_{(2^\kappa)^+}$.   
 
 Our paper is an extension and generalization of the above results of Makkai and Shelah, but with assumptions closer to those of Grossberg and VanDieren.  We introduce a notion that, like the one from \cite{makkaishelah}, is an analogue of the first order notion of coheir.  One of our main results is that, given certain model theoretic assumptions, this notion is in fact an independence notion.

\begin{theorem*}[\ref{fulltheorem}] 
Let $K$ be an AEC with amalgamation, joint embedding, and no maximal models.  If there is some $\kappa > LS(K)$ so that
\begin{enumerate}

	\item $K$ is fully $< \kappa$-tame;
	
	\item $K$ is fully $< \kappa$-type short;
	
	\item $K$ doesn't have an order property; and 
	
	\item $\dnf$ satisfies existence and extension,

\end{enumerate}

then $\dnf$ is an independence relation.
\end{theorem*}

Sections \ref{prelimsection}, \ref{axiomssection}, and \ref{connectingsection} give precise definitions and discussions of the terms in the theorem.  After proving the main theorem in Section \ref{fullsitsection}, we delve into a discussion of the assumptions used; the difficulties of the proof compared with first-order; and some examples of where the theorem can be applied.  At the referee's suggestion, we have included a list of the other main results in the paper; these all take place under the hypothesis of amalgamation, joint embedding, and no maximal models and $\kappa > LS(K)$.

\begin{theorem*}
\begin{itemize}
	\item If $K$ has no weak $\kappa$-order property, is categorical in $\lambda \geq \kappa$, and Galois stable in $\kappa$, then $\kappa^*_\omega(\dnf) = \omega$; see Theorem \ref{ulcthm}.
	\item If $K$ is fully $<\kappa$-tame and -type short, doesn't have the weak $\kappa$-order property, is categorical in $\lambda > \kappa$, and $\dnf$ satisfies existence and extension, then $K_{[\kappa, \lambda)}$ has unique limit models in each cardinality; see Corollary \ref{ulmcor}.
	\item If $\adnf$ is an independence relation, then the forking rank characterizes nonforking on ordinal-ranked Galois types; see Theorem \ref{nf=U-rank}.
	\item If $\kappa$ is strongly compact and $K$ is categorical in $\lambda = \lambda^{<\kappa}$, then $\dnf$ is an independence relation with superstable-like local character; see Theorem \ref{sceasy}.(3).
\end{itemize}
\end{theorem*}

    The main theorem also generalizes the work of \cite{measure2}.  This has several advantages over previous work.
\begin{enumerate}

	\item We generalize Makkai and Shelah \cite{makkaishelah} in two key ways.  First, we replace $L_{\kappa, \omega}$ classes and syntactic types by AECs and Galois types.  This allows the results to apply to classes not axiomatizable in $L_{\kappa, \omega}$, such as quasiminimal excellent classes (which require the `Q' quantifier) and the example of Baldwin, Eklof, and Trlifaj \cite{nperp}.  Second, we replace large cardinal axioms by purely model theoretic hypotheses: tameness and type shortness.  Section \ref{fullsitsection} gives some ZFC examples of AECs with these properties.  Together, tameness and type shortness give a locality condition for when an injection with domain not necessarily a model can be extended into a $K$ embedding;  see \cite{tamelc}.\S3 for a longer discussion.

	\item When reduced to the special case where K has a first-order axiomatization by a complete theory, the independence relation we introduce is coheir, which is equivalent to first-order forking in the stable case we consider.  This is unlike the relatives of splitting and strong splitting used by \cite{hyttinenkesala2006} and \cite{GrLe1}.  This allows us to mimic some first order arguments; see Section 6.
	
	\item  Motivated by a test question of Grossberg\footnote{This first appeared in Grossberg's 1981 MSc thesis.}, Shelah \cite{sh576}, \cite{sh600}, and \cite{Sh705} and Jarden and Shelah \cite{jrsh875} have dealt with the problem whether $I(\lambda,K)=1\leq I(\lambda^+,K)<2^{\lambda^+}$ implies existence of model of cardinality $\lambda^{++}$.  While this question is still open (even under strong set-theoretic assumptions), Shelah managed to get several approximations.  For this, he needed to discover and develop a very rich conceptual infrastructure that occupies more than 500 pages.  One of the more important notions is that of good $\lambda$-frame.  This is a forking-like relation defined using Galois-types over models of cardinality $\lambda$.  Our approach is orthogonal to Shelah's recent work on good $\lambda$-frames and we manage to obtain a forking notion on the class of all models above a natural threshold size, instead of models of a single cardinality. Instead of using $I(\lambda^{+n},K)<2^{\lambda^{+n}}$ for all $1\leq n<\omega$, we assume the lack of an order property, which follows from few models in a single big cardinal.  Unlike Shelah, our treatment does not make use of diamond-like principles as we work in ZFC.  More comparison is given after Definition \ref{indepreldef}.
		
\end{enumerate}

The main theorem allows us to shed light on other questions in the classification theory of AECs, especially concerning superstability.  Superstability in AECs suffers from ``schizophrenia'' \cite[p. 19]{shelahaecbook}.  The main result shows a connection between two conditions, namely that a local character $\kappa^*(\dnf) = \omega$ implies uniqueness of limit models; see Section \ref{lcsection}.  Another question involves properties of ranks (see \cite[Question 3]{Vi}) and Section \ref{uranksection} introduces a rank based on the properties of the independence relation that is ordinal valued exactly when a local character for $\dnf$ holds.

%In a meeting at AIM that was dedicated to Classification Theory for AECs, there was a problem session moderated by Andres Villaveces \cite{Vi}. John Baldwin asked ``Does Shelah's rank satisfy the Lascar inequalities, or is there another rank which does?'' (in the context of Shelah's excellent classes).  Theorem \ref{nf=U-rank} provides an affirmative answer (for a much wider context).  Another question asked by Baldwin and Grossberg at that meeting was ``What is superstability for AECs?''.  While several approximations were offered by various authors, Theorem \ref{U rank bounded iff nf wf} provide the best known approximation to this question.

 Unfortunately, there is no free lunch and we pay for this luxury.  Our payment is essentially in assuming tameness and type-shortness.  As was shown by Boney in \cite{tamelc}, these assumptions are corollaries of certain large cardinal axioms, including the one assumed by Makkai and Shelah; indeed, more recent work of Boney and Unger \cite{bonung} show that assuming \emph{every} AEC is tame is equivalent to a large cardinal axiom.  However, many natural AECs are tame, and it seems to be plausible that tameness and type shortness will be derived in the future from categoricity above a certain Hanf number that depends only on $LS(K)$.

After circulating early drafts of this paper, some of our results were used by Sebastien Vasey \cite{vaseybehemoth} and \cite{vaseyc}, e. g., Theorem \ref{ulcthm} is crucial in his proof of superstability from categoricity.  Also Theorem \ref{sceasy}.(2) is used in a forthcoming article by Grossberg and Vasey \cite{gvas}.

Section 2 gives the necessary background information for AECs.  Section 3 gives a list of common axioms for independence relations and defines the forking relation that we will consider in this paper.  Section 4 gives a fine analysis of when parameterized versions of the axioms from Section 3 hold about our forking relation.  Section 5 gives the global assumptions that make our forking relation an independence relation.  Section 6 introduces a notion that generalizes heir and deduces local character of our forking from this and categoricity.  Section 7 introduces a $U$ rank and shows that it is well behaved.  Section 8 continues the study of large cardinals from \cite{tamelc} and shows that large cardinal assumptions simplify many of the previous sections.

This paper was written while the first author was working on a Ph.D. under the direction of Rami Grossberg at Carnegie Mellon University, and he would like to thank Professor Grossberg for his guidance and assistance in his research in general and relating to this work specifically.  He would also like to thank his wife Emily Boney for her support.  A preliminary version of this paper was presented in a seminar at Carnegie Mellon and we appreciate comments from the participants, in particular Jose Iovino.  We would also like to thank John Baldwin, Adi Jarden, Sebastien Vasey and Andres Villaveces. We are very grateful to the referee 
for comments on several versions of this paper, the referee's reports significantly improved  our paper.  
%%%%%%%%%%%%%%%%%%%%%%%%%%%%%%%%%%%%%%%%%%%%%%%
\section{Preliminaries} \label{prelimsection}
%%%%%%%%%%%%%%%%%%%%%%%%%%%%%%%%%%%%%%%%%%%%%%%

The definition of an Abstract Elementary Class was first given by Shelah in \cite{sh88}.  The definitions and concepts in this section are all part of the literature; in particular, see the books by Baldwin \cite{baldwinbook} and Shelah \cite{shelahaecbook}, the  article by Grossberg \cite{grossberg2002}, or the forthcoming book by Grossberg \cite{ramibook} for general information.

\begin{defin}
We say that $(K, \prec_K)$ is an Abstract Elementary Class iff
\begin{enumerate}

    \item There is some language $L = L(K)$ so that every element of $K$ is an $L$-structure;   

    \item $\prec_K$ is a partial order on $K$;

    \item for every $M, N \in K$, if $M \prec_K N$, then $M$ is an $L$-substructure of $ N$;

    \item $(K, \prec_K)$ respects $L$ isomorphisms: if $f: N \to N'$ is an $L$ isomorphism and $N \in K$, then $N' \in K$ and if we also have $M \in K$ with $M \prec_K N$, then $f(M) \prec_K N'$;

    \item \emph{(Coherence)} if $M_0, M_1, M_2 \in K$ with $M_0 \prec_K M_2$, $M_1 \prec_K M_2$, and $M_0 \subseteq M_1$, then $M_0 \prec M_1$;

    \item \emph{(Tarski-Vaught axioms)} suppose $\seq{M_i \in K : i < \alpha}$ is a $\prec_K$-increasing continuous chain, then

        \begin{enumerate}

            \item $\cup_{i < \alpha} M_i \in K$ and, for all $i < \alpha$, we have $M_i \prec_K \cup_{i < \alpha} M_i$; and

            \item if there is some $N \in K$ so that, for all $i < \alpha$, we have $M_i \prec_K N$, then we also have $\cup_{i < \alpha} M_i \prec_K N$.

        \end{enumerate}

    \item \emph{(Lowenheim-Skolem number)} $LS(K)$ is the first infinite cardinal $\lambda \geq |L(K)|$ such that for any $M \in K$ and $A \subset M$, there is some $N \prec_K M$ such that $A \subset N$ and $\|N\| \leq |A| + \lambda$.

\end{enumerate}

\end{defin}

\begin{remark}
 As is typical, we drop the subscript on $\prec_K$ when it is clear from context and abuse notation by calling $K$ an AEC when we mean $(K, \prec_K)$ is an AEC.  Also, we follow the convention of Shelah that, for $M \in K$, we denote the cardinality of its universe by $\|M\|$.  Also, in this paper, $K$ is always an AEC that has no models of size smaller than the Lowenheim-Skolem number.
\end{remark}

We will briefly summarize some of the necessary basic notations, definitions, and results for AECs; the references contain a more detailed description and development.

\begin{defin}
\begin{enumerate}

	\item[]

	\item A \emph{$K$ embedding} from $M$ to $N$ is an injective $L(K)$-morphism $f: M \to N$ so $f(M) \prec_K N$.

	\item \begin{eqnarray*}
	K_\lambda &=& \{ M \in K : \|M\| = \lambda \}\\
	K_{\leq \lambda} &=& \{ M \in K : \|M\| \leq \lambda \}
	\end{eqnarray*}

	\item $K$ has the \emph{amalgamation property} (AP) iff for any $M \prec N_0, N_1 \in K$, there is some $N^* \in K$ and $f_i: M \to N_i$ such that
\[
 \xymatrix{\ar @{} [dr] N_1  \ar[r]^{f_1}  & N^*\\
M \ar[u] \ar[r] & N_0 \ar[u]_{f_2}
 }
\]
commutes.
	
	\item $K$ has the \emph{joint embedding property} (JEP) iff for any $M_0, M_1 \in K$, there is some $M^* \in K$ and $f_i: M_i \to M^*$.

	\item $K$ has \emph{no maximal models} iff for all $M \in K$, there is $N \in K$ so $M \precneqq N$.

\end{enumerate}
\end{defin}

Note that it is a simple exercise to show that,  if $K$ has the joint embedding property, then having arbitrarily large models is equivalent to having no maximal models.
%see Baldwin, Larson, and Shelah \cite{bls1003}.

We will use the above three assumptions in tandem throughout this paper.  This allows us to make use of a monster model, as in the complete, first order setting; \cite[Section 4.4]{ramibook} gives details.  The monster model $\sea$ is of large size and is universal and model homogeneous for all models that we consider.  As is typical, we assume all elements come from the monster model.  

We use a monster model to streamline our treatment.  However, amalgamation is the only one of the properties that is crucial because it simplifies Galois types.  Joint embedding and no maximal models are rarely used; one major exception is Proposition \ref{orderprop} in the discussion of the order property.  After giving the definition of nonforking in the next section, we briefly detail the differences when we are not working in the context of a monster model.

In AECs, a consistent set of formulas is not a strong enough definition of type; any of the examples of non-tameness will be an example of this and it is made explicit in \cite{untame}.  However, Shelah isolated a semantic notion of type in \cite{sh300} that Grossberg suggested to call \emph{Galois type} in \cite{grossberg2002}.  In his book \cite{shelahaecbook}, Shelah calls this \emph{orbital type}. This notion replaces the first order notion of type for AECs.  Crucially, we allow Galois types of infinite lengths.

\begin{defin} Let $K$ be an AEC, $\lambda \geq LS(K)$, and $I$ be a nonempty set.
 \begin{enumerate}
  
	\item Let $M$ in $K$ and $\seq{a_i \in \sea : i \in I}$ be a sequence of elements.  The Galois type of $\seq{a_i \in \sea : i \in I}$ over $M$ is denoted $gtp(\seq{a_i \in \sea : i \in I}/M)$ and is the orbit of $\seq{a_i :i \in I}$ under the action of automorphisms of $\sea$ fixing $M$.  That is, $\seq{a_i \in \sea : i \in I}$ and $\seq{b_i : i \in I}$ have the same Galois type over $M$ iff there is $f \in Aut_M \sea$ so that $f(a_i) = b_i$ for all $i \in I$.
    
     \item For $M \in K$, $gS^I(M) = \{ gtp(\seq{a_i : i \in I}/M) : a_i \in \sea $ for all $i \in I \}$.
    
    \item Suppose $p = gtp(\seq{a_i : i \in I}/M)  \in gS^I(M)$ and $N \prec M$ and $I' \subset I$.  Then, $p \rest N \in gS^I(N)$ is $gtp(\seq{a_i : i \in I}/N)$ and $p^{I'} \in gS^{I'}(M)$ is $ gtp(\seq{a_i : i \in I'}/M)$.

	\item Given a Galois type $p \in gS^I(M)$, then the domain of $p$ is $M$ and the length of $p$ is $I$.
	
	\item If $p = gtp(A/M)$ is a Galois type and $f \in Aut $ $\sea$, then $f(p) = gtp(f(A)/f(M))$.

 \end{enumerate}

\end{defin}

\begin{remark}
We sometimes write that the type of two sets (say $X$ and $Y$) are equal; given the above definitions, this really means there is some enumeration $X = \seq{x_i : i \in I}$ and $Y = \seq{y_i : i \in I}$ so that the types of the sequences are equal.  If we reference some $x \subset X$ and the `corresponding part of $Y$,' then this refers to $y \subset Y$ indexed by the same set that indexes $x$; that is, $y = \seq{y_i : i \in I \te{ and } x_i \in x}$.
\end{remark}

Along with types comes a notion of saturation, called Galois saturation.  A degree of Galois saturation will be necessary when we deal with our independence relation, so we offer a definition here.  Additionally, we include a lemma of Shelah that characterizes saturation by model homogeneity.

\begin{defin}
\begin{enumerate}

	\item A model $M \in K$ is $\mu$-Galois saturated iff for all $N \prec M$ such that $\|N\| < \mu$ and $p \in gS(N)$, we have that $p$ is realized in $M$.
	
	\item A model $M \in K$ is $\mu$-model homogeneous iff for all $N \prec M$ and $N' \succ N$ such that $\|N'\| < \mu$, there is $f:N' \to_N M$.

\end{enumerate}
\end{defin}

\begin{lemma}[\cite{sh576}.0.26.1]
Let $\lambda >  LS(K)$ and $M \in K$ and suppose that $K$ has the amalgamation property.  Then $M$ is $\lambda$-Galois saturated iff $M$ is $\lambda$-model homogeneous.
\end{lemma}

We conclude the preliminaries by recalling two locality properties that are key for this paper: tameness and type shortness.  Tameness was first isolated by Grossberg and VanDieren \cite[Definition 3.2]{tamenessone}, although a weaker version had been used by Shelah \cite{sh394} in the midst of a proof.  Later, Grossberg and VanDieren \cite{tamenesstwo} \cite{tamenessthree} showed that a strong form of Shelah's Categoricity Conjecture holds for tame AECs.  Type shortness was defined by the first author in \cite[Definition 3.3]{tamelc} as a dual property for tameness.  There, type shortness and tameness were derived from large cardinal hypotheses.  \\
Recall that $P_\kappa I = \{ I_0 \subset I: |I_0| < \kappa\}$ and $P^*_\kappa M = \{ M^- \prec M : \|M^-\| < \kappa \}$.

\begin{defin}\label{tame-def}
\begin{enumerate}
	\item $K$ is \emph{$(<\kappa, \lambda)$-tame for $\theta$-length types} iff, for all $p, q \in gS^I(M)$ with $\|M\| = \lambda$ and $|I| = \theta$, we have $p = q$ iff, for all $M^- \in P^*_\kappa M$, $p \rest M^- = q \rest M^-$.

	\item $K$ is \emph{$(<\kappa, \lambda)$-type short for $\theta$-sized domains} iff, for all $p, q \in gS^I(M)$ with $|I| = \lambda$ and $\|M\| = \theta$, we have $p = q$ iff, for all $I_0 \in P_\kappa I$, $p^{I_0} = q ^{I_0}$.

	\item K is \emph{fully $<\kappa$-tame and -type short} iff it is $(<\kappa,\lambda)$-tame for $\theta$-length types and $(<\kappa, \lambda)$-type short for $\theta$-sized domains for all $\theta$ and all $\lambda \geq \kappa$.
\end{enumerate}
\end{defin}

We parameterize the properties to get exact results in Proposition \ref{lvlanalprop}, but the reader can focus on the ``fully $<\kappa$-tame and -type short'' case if desired.  Note, by \cite[Theorem 3.5]{tamelc}, the type shortness implies the tameness.  However, we include both hypotheses for clarity.

%%%%%%%%%%%%%%%%%%%%%%%%%%%%%%%%%%%%%%%%%%%%%%%
\section{Axioms of an independence relation and the definition of forking} \label{axiomssection}
%%%%%%%%%%%%%%%%%%%%%%%%%%%%%%%%%%%%%%%%%%%%%%%

The following hypothesis and definition of non-forking is central to this paper:

\begin{hypothesis}\label{hyp}
Assume that $K$ has no maximal models and satisfies the $\lambda$-joint embedding and $\lambda$-amalgamation properties for all $\lambda \geq LS(K)$.

Fix a cardinal $\kappa > LS(K)$.  The nonforking is defined in terms of this $\kappa$ and all subsequent uses of $\kappa$ will refer to this fixed cardinal, until Section \ref{lcrevisitedsection}.  If we refer to a model, tuple, or Galois type as `small,' then we mean its size is $< \kappa$, its length is of size $< \kappa$, or both its domain and its length are small.
\end{hypothesis}

\begin{defin}
Let $M_0 \prec N$ be models and $A$ be a set.  We say that \emph{$gtp(A/N)$ does not fork over $M_0$}, written $A \dnf_{M_0} N$, iff for all small $a \in A$ and all small $N^- \prec N$, we have that $tp(a/N^-)$ is realized in $M_0$.  We call this $\kappa$-coheir or $<\kappa$-satisfiability. 
\end{defin}

That is, a type does not fork over a base model iff all small approximations to it, both in length and domain, are realized in the base model.  This definition is a relative of the finite satisfiability--also known as coheir--notion of forking that is extensively studied in stable theories.  It is an AEC version of the non-forking defined in Makkai and Shelah \cite[Definition 4.5]{makkaishelah} for categorical $L_{\kappa, \omega}$ theories when $\kappa$ is strongly compact.

We now list the properties that, our nonforking notion will have.  These properties can be thought of as axiomatizing an independence relation.  The ones listed below are commonly considered and are similar to the properties that characterize nonforking in first order, stable theories, although this list is most inspired by \cite[Proposition 4.4]{makkaishelah}.  However, many of these properties have been changed because we require the bottom and right inputs to be models.  This is similar to good $\lambda$-frames, which appear in \cite[Section II.2]{shelahaecbook}, although we don't require the parameter set $A$ to be a singleton and we allow the sets and models to be of any size.

The properties we introduce are heavily parameterized.  The interesting and hard to prove properties--Existence, Extension, Uniqueness, and Symmetry--are each given with parameters among $\lambda$, $\mu$, $\chi$, and $\theta$.  These parameters allows us to conduct a fine analysis of exactly what assumptions are required to derive these properties.  The order of these parameters is designed to be as uniformized as possible: if appropriate, when referring to $A \dnf_{M_0} N$,  the $\lambda$ refers to the size of $A$, $\mu$ refers to the size of $M_0$, and $\chi$ refers to the size of $N$.  If we write a property without parameters, then we mean that property for all possible parameters.

\begin{defin} \label{indepreldef}
Fix an AEC $K$.  Let $\adnf$ be a ternary relation on models and sets so that $A \adnf_{M_0} N$ implies that $A$ is a subset of the monster model and $M_0 \prec N$ are both models.  We say that $\adnf$ is an \emph{independence relation} iff it satisfies all of the following properties for all cardinals referring to sets and all cardinals that are at least $\kappa$ when the cardinal refers to a model.
\begin{enumerate}

		\item[($I$)] {\bf Invariance}\\
		Let $f \in Aut$ $\sea$ be an isomorphism.  Then $A \adnf_{M_0} N$ implies $f(A) \adnf_{f(M_0)} f(N)$.
		
		\item[$(M)$] {\bf Monotonicity}\\
		If $A \adnf_{M_0} N$ and $A' \subset A$ and $M_0 \prec M_0' \prec N' \prec N$, then $A' \adnf_{M_0'} N'$.
			
		\item[$(T)$]  {\bf Transitivity}\\
		If $A \adnf_{M_0'} N$ and ${M_0'} \adnf_{M_0} N$ with $M_0 \prec M_0'$, then $A \adnf_{M_0} N$.
		
		\item[$(C)_{< \kappa}$] {\bf Continuity}
		\begin{enumerate}
			
			\item If for all small $A' \subset A$ and small $N' \prec N$, there are $M_0' \prec M_0$ and $N' \prec N^* \prec N$ such that $M_0' \prec N^*$ and  $A' \adnf_{M_0'} N^*$, then $A \adnf_{M_0} N$.
			
			\item If $I$ is a $\kappa$-directed partial order and $\seq{A_i, M_0^i \mid i \in I}$ are increasing such that $A_i \adnf_{M_0^i} N$ for all $i \in I$, then $\cup_{i \in I}A_i \adnf_{\cup_{i \in I} M^i_0} N$.

		\end{enumerate}

		\item[$(E)_{(\lambda, \mu, \chi, \theta)}$] 
		
	\begin{enumerate}
	
	\item  
	{\bf Existence}\\
		Let $A$ be a set and $M_0$ be a model of sizes $\lambda$ and $\mu$, respectively.  Then  $A \adnf_{M_0} M_0$.
		\item
		{\bf Extension}\\
		Let $A$ be a set and $M_0$ and $N$ be models of sizes $\lambda$, $\mu$, and $\chi$, respectively, so that $M_0 \prec N$ and $A \adnf_{M_0} N$.  If $N^+ \succ N$ of size $\theta$, then there is $A'$ so $A' \adnf_{M_0} N^+$ and $gtp(A'/N) = gtp(A/N)$.

		\end{enumerate}

		\item[$(S)_{(\lambda, \mu, \chi)}$] {\bf Symmetry}\\
		Let $A_1$ be a set, $M_0$ be a model, and $A_2$ be a set of sizes $\lambda$, $\mu$, and $\chi$, respectively, so that there is a model $M_2$ with $M_0 \prec M_2$ and $A_2 \subset M_2$ such that $A_1 \adnf_{M_0} M_2$.  Then there is a model $M_1 \succ M_0$ that contains $A_1$ such that $A_2 \adnf_{M_0} M_1$.
		
		\item[$(U)_{(\lambda, \mu, \chi)}$] {\bf Uniqueness}\\
		Let $A$ and $A'$ be sets and $M_0 \prec N$ be models of sizes $\lambda$, $\lambda$, $\mu$, and $\chi$, respectively.  If $gtp(A/M_0) = gtp(A'/M_0)$ and $A \adnf_{M_0} N$ and $A' \adnf_{M_0} N$, then $tp(A/N) = tp(A'/N)$.

\end{enumerate}

\end{defin}

A discussion of these axioms and their relation to other nonforking notions is in order.

As mentioned in the introduction, this notion is somewhat orthogonal to Shelah's notion of good $\lambda$-frames (see \cite[Definition II.2.1]{shelahaecbook}).  While both attempt to axiomatize a nonforking relation, we allow greater generality by considering Galois types of arbitrary length over all models of a sufficiently large size.  In contrast, Shelah deals with a subclass of unary Galois types only over a fixed size $\lambda$.  On the other hand, since good $\lambda$-frames attempt to axiomatize superstability (rather than stability as we do here), good $\lambda$-frames have stronger continuity and local character properties.

Notice that the Existence property implies that $M_0\in K_{\geq \kappa}$ is $\kappa$-Galois saturated when $\dnf$ is $\kappa$-coheir.  However, this is not a serious restriction in comparison with other work.  In other cases where AECs have a strong nonforking notions, some level of categoricity is assumed; see Makkai and Shelah \cite{makkaishelah}, Shelah \cite[Section II.3]{shelahaecbook}, and Vasey \cite{vaseytameframe}.  The categoricity hypotheses are used to ensure large amounts of Galois saturation, as we have here.  Indeed, Theorem \ref{catexist} below shows that categoricity in some $\lambda = \lambda^{< \kappa}$ implies that all sufficiently large models are $\kappa$-Galois saturated.

The monotonicity and invariance properties are actually necessary to justify our formulation of nonforking as based on Galois types.  Without them, whether or not a Galois type doesn't fork over a base model could depend on the specific realization of the chosen type.  Since these properties are clearly satisfied by our definition of nonforking, this is not an issue.  

The axiom $(E)_{(\lambda, \mu, \chi, \theta)}$ combines two notions.  The first is Existence: that a Galois type does not fork over its domain.  This is similar to the consequence of simplicity in first order theories that a type does not fork over the algebraic closure of its domain.  As mentioned above, in this context, existence is equivalent to every model being $\kappa$-Galois saturated.  In the first order case, where finite satisfiability is the proper analogue of our non-forking, existence is an easy consequence of the elementary substructure relation.  In \cite{makkaishelah}, this holds for $<\kappa$ satisfiability, their nonforking, because types are formulas from $L_{\kappa, \kappa}$ and, due to categoricity, the strong substructure relation is equivalent to $\prec_{L_{\kappa, \kappa}}$.

The second notion is the extension of nonforking Galois types.  In first order theories (and in \cite{makkaishelah}), this follows from compactness, but is more difficult in a general AEC.  We have separated these notions for clarity and consistency with other sources, but could combine them in the following statement of $(E)_{(\lambda, \mu, \mu, \chi)}$:
\begin{enumerate}
\item[] Let $A$ be a set and $M_0$ and $N$ be models of sizes $\lambda$, $\mu$, and $\chi$, respectively, so that $M_0 \prec N$.  Then there is some $A'$ so that $gtp(A'/M_0) = gtp(A/M_0)$ and $A' \adnf_{M_0} N$.
\end{enumerate}
As an alternative to assuming $(E)$, and thus assuming all models are $\kappa$-Galois saturated, we could simply work with the definition and manipulate the nonforking relationships that occur.  This is the strategy in Section \ref{lcsection}.  In such a situation, $\kappa$-Galois saturated models, which will exist in $\lambda^{<\kappa}$, will satisfy the existence axiom.%\footnote{WB: This was one of the places that I was thinking of ``responding to Sebastien.''  Something along the lines of ``not initially clear how strong.  however, recent results give classes with it'' or ``vasey was able to derive it for certain classes from categoricity''}

%The extension property has a long history.  At first was made by Buechler and Lessmann \cite{BuLe} in the context of homogenous model theory and is called there simplicity, this tradition was continued by several authors including Hyttinen and Kesala \cite{}.  For a several years it was really strange to assume simplicity without having a proof that stability (in homogenous classes) implies simplicity.  Eventually Shelah constructed an example of a simple but not stable class (see the last result of \cite{hyttinenlessman2002}).

%Sebastien Vasey in a tour de force (consisting of about 125 pages, \cite{vaseybehemoth} and \cite{vaseyc}), in part building on this paper managed to conclude a similar result to our main theorem by replacing the extension property by several strong (but natural assumptions). A simplified version of his result  (Theorem 9.3 from \cite{vaseyc}) is: Suppose there is a $\kappa=\beth_\kappa>LS(K)$, $K$ is $<\kappa$-tame and short and $K$ is $\lambda$-stable for some $\lambda\geq \kappa$ then our the $\kappa$-cohier relation is an independence relation on the class $\{M\in K\mid M \text{ is } \lambda^+\text{-staturated} \}$.

The relative complexity of the symmetry property is necessitated by the fact that the right side object is required to be a model that contains the base.  If the left side object already satisfied this, then there is a simpler statement.

\begin{prop}
If $(S)_{(\lambda, \mu, \chi)}$ holds, then so does the following
\begin{enumerate}
	\item[$(S^*)_{(\lambda, \mu, \chi)}$] Let $M$, $M_0$, and $N$ be models of size $\lambda$, $\mu$, and $\chi$, respectively so that $M_0 \prec N$ and $M_0 \prec M$.  Then $M \adnf_{M_0} N$ iff $N \adnf_{M_0} M$.
\end{enumerate}
\end{prop}

In first order stability theory, many of the key dividing lines depend on the local character $\kappa(T)$, which is the smallest cardinal so that any type doesn't fork over some subset of its of domain of size less than $\kappa(T)$.  The value of this cardinal can be smaller than the size of the theory, e.g. in an uncountable, superstable theory.  However, since types and nonforking occur only over models, the smallest value the corresponding cardinal could take would be $LS(K)^+$.  This is too coarse for many situations, although see Boney and Vasey \cite{unionsat} for a comparison.  Instead, we follow \cite{shvi635}, \cite[Chapter II]{shelahaecbook}, and \cite{tamenessone} by defining a local character cardinal based on the length of a resolution of the base rather than the size of cardinals.  As different requirements appear in different places, we give two definitions of local character: one with no additional requirement, as in \cite[Chapter II]{shelahaecbook}, and one  requiring that successor models be universal, as in \cite{shvi635} and \cite{tamenessone}.

\begin{defin} \label{lc-def}
$\kappa_\alpha(\adnf) = \min \{ \lambda \in REG \cup \{\infty\} :$ for all regular $\mu \geq \lambda$ and all increasing, continuous chains $ \seq{M_i : i < \mu}$  and all sets $A$ of size less than $\alpha$, there is some $ i_0 < \mu $ so $ A \adnf_{M_{i_0}} \cup_{i < \mu} M_i \}$

$\kappa^*_\alpha(\adnf) = \min \{ \lambda \in REG \cup \{\infty\} :$ for all $\mu = \cf \mu \geq \lambda$ and all increasing, continuous chains$ \seq{M_i : i < \mu}$ with $ M_{i+1} $ universal over $ M_i $ which is $\kappa$-saturated and all sets $A$ of size less than $\alpha$, there is some $ i_0 < \mu $ so $ A \adnf_{M_{i_0}} \cup_{i < \mu} M_i \}$

In either case, if we omit $\alpha$, then we mean $\alpha = \omega$.
\end{defin}

In Section \ref{lcsection}, we return to these properties and find a natural, sufficient condition that implies that $\kappa^*(\dnf) = \omega$.

Although we do not do so here, the notions of $\kappa$-coheir makes sense in AECs with some level of amalgamation but without the full strength of a monster model.  In such an AEC,  the definition of the Galois type of $A$ over $N$ must be augmented by a model containing both; that is, some $\bigM \in K$ so $A \subset \bigM$ and $N \prec \bigM$.  We denote this type $gtp(A/N, \bigM)$.  Similarly, we must add this fourth input to the nonforking relation that contains all other parameters.  Then $A \dnf^{\bigM}_{M_0} N$ iff $M_0 \prec N \prec \bigM$ and $A \subset |\bigM|$ and all of the small approximations to the Galois type of $A$ over $N$ \emph{as computed in $\bigM$} are realised in $M_0$.  The properties are expanded similarly with added monotonicity for changing the ambient model $\bigM$ and the allowance that new models that are found by properties such as existence or symmetry might exist in a larger big model $\bigN$.  All theorems proved in this paper about nonforking only require amalgamation, although some of the results referenced make use of the full power of the monster model.

We end this section with an easy exercise in the definition of nonforking that says that $A$ and $N$ must be disjoint outside of $M_0$.

\begin{prop}
If we have $A \dnf_{M_0} N$, then $A \cap |N| \subset |M_0|$.
\end{prop}

{\bf Proof:} Let $a \in A \cap |N|$.  Since $N$ is a model, we can find a small $N^- \prec N$ that contains $a$.  Then, by the definition of nonforking, $gtp(a/N^-)$ must be realized in $M_0$.  But since $a \in |N^-|$, this type is algebraic so the only thing that can realize it is $a$.  Thus, $a \in |M_0|$. \hfill \dag

%%%%%%%%%%%%%%%%%%%%%%%%%%%%%%%%%%%%%%%%%%%%%%%
\section{Connecting Existence, Symmetry and Uniqueness} \label{connectingsection}
%%%%%%%%%%%%%%%%%%%%%%%%%%%%%%%%%%%%%%%%%%%%%%%

In this section, we investigate what AEC properties cause the axioms of our independence relation to hold; recall that we are working under Hypothesis \ref{hyp} that $K$ is an AEC with amalgamation, joint embedding, and no maximal models and that $\kappa > LS(K)$ is fixed.  The relations are summarized in the proposition below.

\begin{prop} \label{lvlanalprop}
Suppose that $K$ doesn't have the   weak $\kappa$-order property and is $(< \kappa, \lambda + \chi)$-type short for $\theta$-sized domains and $(< \kappa, \theta)$-tame for $<\kappa$ length types.  Then, for the $\kappa$-coheir relation,
\begin{enumerate}

	\item $(E)_{(\chi, \theta, \theta, \lambda)}$ implies $(S)_{(\lambda, \theta, \chi)}$.

	\item $(S)_{(<\kappa, \theta, <\kappa)}$ implies $(U)_{(\lambda, \theta, \chi)}$.
		
\end{enumerate}
\end{prop}
Recall Definition \ref{tame-def} for tameness and type shortness.  This proposition and the lemma used to prove it below rely on an order property.

\begin{defin}\label{orderdef}
$K$ has the \emph{weak $\kappa$-order property} iff there are lengths $\alpha, \beta < \kappa$, a model $M \in K_{< \kappa}$, and types $p \neq q \in gS^{\alpha + \beta}(M)$ such that there are sequences $\seq{a_i \in {}^\alpha \sea : i < \kappa}$ and $\seq{b_i \in {}^\beta \sea : i < \kappa}$ such that, for all $i, j < \kappa$,
\begin{eqnarray*}
i \leq j \implies gtp(a_i b_j/M) = p\\
i > j \implies gtp(a_i b_j/M) = q
\end{eqnarray*}
\end{defin}

This order property is a generalization of the first order version to our context of Galois types and infinite sequences.  This is one of many order properties proposed for the AEC context and is similar to 1-instability that is studied by Shelah in \cite{sh1019} in the context of $L_{\theta, \theta}$ theories where $\theta$ is strongly compact.  The adjective `weak' is in comparison to the $(<\kappa, \kappa)$-order property in Shelah \cite[Definition 4.3]{sh394}.  The key difference is that \cite{sh394} requires the existence of ordered sequences of \emph{any} length (i.e. the existence of $\seq{a_i, b_i: i <\delta}$ for all ordinals $\delta$), while we only require a sequence of length $\kappa$.  We discuss the implications of the weak $\kappa$-order property in the next section.  For now, we use it's negation to prove the following result, similar to one in \cite[Proposition 4.6]{makkaishelah}, based on first order versions due to Poizat and Lascar.

\begin{lemma}\label{usefullemma}
Suppose $K$ is an AEC that is $(< \kappa, \theta)$-tame for $<\kappa$ length types and doesn't have the weak $\kappa$-order property.  Let $M_0 \prec M, N$ such that $\|M_0\| = \theta$ and let $a, b, a' \in {}^{<\kappa}\sea$ such that $b \in N$ and $a' \in M$.  If
$$gtp(a/M_0) = gtp(a'/M_0) \te{ and } a \dnf_{M_0} N \te{ and } b \dnf_{M_0} M$$
then $gtp(ab/M_0) = gtp(a'b/M_0)$.
\end{lemma}

{\bf Proof:}  Assume for contradiction that $gtp(ab/M_0) \neq gtp(a'b/M_0)$.  We will build sequences that witness the weak $\kappa$-order property.  By tameness, there is some $M^- \prec M_0$ of size $< \kappa$ such that $gtp(ab/M^-) \neq gtp(a'b/M^-)$.  Set $p = gtp(ab/M^-)$ and $q = gtp(a'b/M^-)$.  We will construct two sequences $\seq{a_i \in {}^{\ell(a)}M_0 : i < \kappa}$ and $\seq{b_i \in {}^{\ell(b)} M_0 : i < \kappa}$ by induction.  We will have, for all $i < \kappa$
\begin{enumerate}

	\item $a_i b \vDash p$;
	\item $a_i b_j \vDash q$ for all $j < i$;
	\item $a b_i \vDash q$; and
	\item $a_i b_j \vDash p$ for all $j \geq i$.

\end{enumerate}
Note that, since $b_i \in {}^{\ell(b)}M_0$, $(3)$ is equivalent to $a' b_i \vDash q$.

{\bf This is enough:} $(2)$ and $(4)$ are the properties necessary to witness the weak $\kappa$-order property.

{\bf Construction:}  Let $i < \kappa$ and suppose that we have constructed the sequences for all $j < i$.  Let $N^+ \prec N$ of size $< \kappa$ contain $b$, $M^-$, and $\{ b_j : j < i\}$.  Because $a \dnf_{M_0} N$, there is some $a_i \in M_0$ that realizes $gtp(a'/N^+)$.  This is witnessed by $f \in Aut_{N^+} \sea$ with $f(a) = a_i$.\\
{\bf Claim:} $(1)$ and $(2)$ hold.\\
$f$ fixes $M^-$ and $b$, so it witnesses that $gtp(ab/M^-) = gtp(a_i b/M^-)$.  Similarly, it fixes $b_j$ for $j < i$, so it witnesses $q = gtp(ab_j/M^-) = gtp(a_i b_j/M^-)$.\hfill $\dag_{Claim}$
Similarly, pick $M^+ \prec M$ of size $< \kappa$ to contain $M^-$, $a'$, and $\{a_j : j \leq i \}$.  Because $b \dnf_{M_0} M'$, there is $b_i \in M_0$ that realizes $gtp(b/M^+)$.  As above, $(3)$ and $(4)$ hold. \hfill \dag\\

Now we are ready to prove our theorems regarding when the properties of $\dnf$ hold.  The first four properties always hold from the definition of nonforking.

\begin{prop} \label{easypropthm}
$\dnf$ satisfies $(I)$, $(M)$, $(T)$, and $(C)_{< \kappa}$.
\end{prop}

To get the other properties, we have to rely on some degree of tameness, type shortness, no weak order property, and the property $(E)$.\\

{\bf Proof of Proposition \ref{lvlanalprop}:} \begin{enumerate}

\item Suppose $(E)_{(\chi, \theta, \theta, \lambda)}$ holds.  Let $A_2 \dnf_{M_0} M_1$ and $A_1 \subset |M_1|$ with $|A_2| = \lambda$, $\|M_0\|=\theta$, and $|A_1| = \chi$; WLOG $\|M_1\| = \chi$.  Let $M_2$ contain $A_2$ and $M_0$ be of size $\lambda$.  By $(E)_{(\chi, \theta, \theta, \lambda)}$, there is some $A_1'$ such that $gtp(A_1/M_0) = gtp(A_1'/M_0)$ and $A_1' \dnf_{M_0} M_2$.  It will be enough to show that $gtp(A_1 A_2/M_0) = gtp(A_1' A_2/M_0)$.  By $(< \kappa, \lambda+\chi)$-type shortness over $\theta$-sized domains, it is enough to show that, for all $a_2 \in A_2$ and corresponding $a_1 \in A_1$ and $a'_1 \in A_1'$ of length $< \kappa$, we have $gtp(a_1 a_2/M_0) = gtp(a_1' a_2/M_0)$.  By $(M)$, we have that $a_1' \dnf_{M_0} M_2$ and $a_2 \dnf_{M_0} M_1$, so this follows by Lemma \ref{usefullemma} above.

Now that we have shown the type equality, let $f \in Aut_{M_0} \sea$ such that $f(A_1 A_2) = A_1' A_2$.  Applying $f$ to $A_1' \dnf_{M_0} M_2$, we get that $A_1 \dnf_{M_0} f(M_2)$ and $A_2 = f(A_2) \subset f(M_2)$, as desired.

\item Suppose $(S)_{(< \kappa, \theta, < \kappa)}$.  Let $A$ and $A'$ be sets of size $\lambda$ and $M_0 \prec N_0$ of size $\theta$ and $\chi$, respectively, so that

$$gtp(A/M_0) = gtp(A'/M_0) \te{ and } A \dnf_{M_0} N \te{ and } A' \dnf_{M_0} N$$

As above, it is enough to show that $gtp(AN/M_0) = gtp(A'N/M_0)$.  By type shortness, it is enough to show this for every $n \in N$ and corresponding $a \in A$ and $a' \in A'$ of lengths less than $\kappa$.  By $(M)$, we know that $a \dnf_{M_0} N$ and $a' \dnf_{M_0} N$.  By applying $(S)_{(<\kappa, \theta, < \kappa)}$ to the former, there is $N_a \succ M_0$ containing $a$ such that $n \dnf_{M_0} N_a$.  As above, Lemma \ref{usefullemma} gives us the desired conclusion. \hfill \dag\\
		
\end{enumerate} 

%%%%%%%%%%%%%%%%%%%%%%%%%%%%%%%%%%%%%%%%%%%%%%%
\section{The main theorem} \label{fullsitsection}
%%%%%%%%%%%%%%%%%%%%%%%%%%%%%%%%%%%%%%%%%%%%%%%

We now state the ideal conditions under which our nonforking works.  We reiterate Hypothesis \ref{hyp} in the statement of the theorem for clarity.

\begin{theorem} \label{fulltheorem} \label{mainthm}
Let $K$ be an AEC with amalgamation, joint embedding, and no maximal models.  If there is some $\kappa > LS(K)$ such that
\begin{enumerate}

	\item $K$ is fully $< \kappa$-tame and -type short;
		
	\item $K$ doesn't have the weak $\kappa$-order property; and
	
	\item $\dnf$ satisfies $(E)$

\end{enumerate}

then $\dnf$ is an independence relation.
\end{theorem}

{\bf Proof:} First, by Proposition \ref{easypropthm}, $\dnf$ always satisfies $(I)$, $(M)$, $(T)$, and $(C)_{< \kappa}$.  Second, $(E)$ is part of the hypothesis.  Third, by the other parts of the hypothesis, we can use Proposition \ref{lvlanalprop}.  Let $\chi$, $\theta$, and $\lambda$ be cardinals.  We know that $(E)_{(\chi, \theta, \theta, \lambda)}$ holds, so $(S)_{(\lambda, \theta, \chi)}$ holds.  From this, we also know that $(S)_{(<\kappa, \theta, <\kappa)}$ holds.  Thus, $(U)_{(\lambda, \theta, \chi)}$ holds.  So $\dnf$ is an independence relation as in Definition \ref{indepreldef}. \hfill \dag\\

In the following sections, we will assume the hypotheses of the above theorem and use $\dnf$ as an independence relation.  First, we discuss the hypotheses and argue for their naturality.

\underline{``amalgamation, joint embedding, and no maximal models''}\\
	These are a common set of assumptions when working with AECs that appear often in the literature; see \cite{sh394}, \cite{tamenessthree}, and \cite{gvv} for examples.  Readers interested in work on AECs \emph{without} these assumptions are encouraged to see \cite{sh576} or Shelah's work on good $\lambda$-frames in \cite{shelahaecbook} and \cite{jrsh875}.
	
\underline{ ``fully $< \kappa$-tame and -type short''}\\
	As discussed in \cite{tamelc}, these assumptions say that Galois types are equivalent to their small approximations.  Without this equivalence, there is no reason to think that our nonforking, which is defined in terms of small approximations, would say anything useful about an AEC.
	
	On the other hand, we argue that these properties will occur naturally in any setting with a notion of independence or stability theory.  The introduction of \cite{tamenessthree} observes that this happens in all known cases.  Additionally, the following proposition says that the existence of a nonforking-like relation that satisfies stability-like assumptions implies tameness and some stability.
	\begin{prop}\label{tameind-prop}
	If there is a nonforking-like relation $\adnf$ that satisfies $(U)$, $(M)$, and $\kappa_\alpha(\adnf) < \infty$, then $K$ is $(< \mu, \mu)$ tame for less than $\alpha$ length types for all regular $\mu \geq \kappa_\alpha(\adnf)$.
	\end{prop}
	
	{\bf Proof:} Let $p \neq q \in gS^{<\alpha}(M)$ so their restriction to any smaller submodel is equal and let $\seq{M_i \in K_{< \mu} : i < \mu}$ be a resolution of $M$.  By the local character, there are $i_p$ and $i_q$ such that $p$ does not fork over $M_{i_p}$ and $q$ does not fork over $M_{i_q}$.  By $(M)$, both of the types don't fork over $M_{i_p + i_q}$ and, by assumption, $p \rest M_{i_p + i_q} = q \rest M_{i_p + i_q}$.  Thus, by $(U)$, we have $p = q$.  \hfill \dag\\

	The results of \cite[Section 3]{tamelc} allow us to get a similar result for type shortness.
	  
	The arguments of \cite[Proposition 4.14]{makkaishelah} show that this can be used to derive stability-like bounds on the number of Galois types.
			
\underline{``no weak $\kappa$ order property''}\\
	In first order model theory, the order property and its relatives (the tree order property, etc) are well-studied as the non-structure side of dividing lines.  In broader contexts such as ours, much less is known.  Still, there are some results, such as Shelah \cite[Chapter III]{shelahnonstructurebook}, which shows that a strong order property, akin to getting any desired order of a certain size in an EM model, implies many models.  Note that Shelah does not explicitly work inside an AEC, but the proofs and definitions are sufficiently general and syntax free to apply here.

	Ideally, the weak $\kappa$-order property could be shown to imply non-structure for an AEC.  While this is not currently known in general, we have two special cases where many models follows by combinatorial arguments and the work of Shelah.
	
	First, if we suppose that $\kappa$ is inaccessible, then we can use Shelah's work to show that there are the maximum number of models in every size above $\kappa$.  We will show that, given any linear order, there is an EM model with the order property for that order.  This implies \cite{shelahnonstructurebook}'s notion of ``weakly skeleton-like",  which then implies many models by \cite[Conclusion III.3.25]{shelahnonstructurebook}.
	\begin{prop}\label{orderprop}
	Let $\kappa$ be inaccessible and suppose $K$ has the weak $\kappa$-order property.  Then, for all linear orders $I$, there is EM model $M^*$, small $N \prec M^*$, $p \neq q \in gS(M)$, and $\seq{a_i, b_i \in M^* : i \in I}$ such that, for all $i, j \in I$,
	\begin{eqnarray*}
	i \leq j \implies gtp(a_i b_j/M) = p\\
	i > j \implies gtp(a_i b_j/M) = q
	\end{eqnarray*}
	Thus, for all $\chi > \kappa$, $K_\chi$ has $2^\chi$ nonisomorphic models.
	\end{prop}
	We sketch the proof and refer the reader to \cite{shelahnonstructurebook} for more details.\\
	{\bf Proof Outline:}  Let $p \neq q \in gS(N)$ and $\seq{a_i, b_i : i < \kappa}$ witness the weak order property.  Since $K$ has no maximal models, we may assume that this occurs inside an EM model.  In particular, there is some $\Phi$ proper for linear orders so $N \prec EM(\kappa, \Phi) \rest L$ that contains $\seq{a_i, b_i : i < \kappa}$, $L(\Phi)$ contains Skolem functions, and $\kappa$ is indiscernible in $EM(\kappa, \Phi)\rest L$.  Recall that, for $X \subset EM(\kappa, \Phi)$, we have $\cont(X) := \cap \{I \subset \kappa : X \subset |EM(I, \Phi)| \}$.  By inaccessibility, we can thin out $\{ \cont(a_i b_i): i < \kappa\}$ to $\{ \cont(a_i b_i) : i \in J\}$ that is a head-tail $\Delta$ system of size $\kappa$ and are all generated by the same term and have the same quantifier free type in $\kappa$.  Since $\kappa$ is regular and $\cont(N)$ is of size $< \kappa$, we may further assume the non-root portion of this $\Delta$ system is above $\sup (\cont(N))$.
	
	By the definition of EM models, we can put in any linear order into $EM( \cdot, \Phi) \rest L$ and get a model in $K$.  Thus, we can take the blocks that generate each $a_i b_i$ with $i \in J$ and arrange them in any order desired.  In particular, we can arrange them such that they appear in the order given by $I$.  Then, the order indiscernibility implies that the order property holds as desired.
	
	We have shown the hypothesis of \cite[Conclusion III.3.25]{shelahnonstructurebook} and the final part of our hypothesis is that theorem's conclusion. \hfill \dag\\
	
	We can also make use of these results without large cardinals.  To do so, we `forget' some of the tameness and type shortness our class has to get a slightly weaker relation.  Suppose $K$ is $< \kappa'$ tame and type short.  Let $\lambda$ be regular such that $\lambda^{\kappa'} = \lambda > \kappa'$.  By the definitions, $K$ is also $<\lambda$ tame and type short, so take $\lambda$ to be our fixed cardinal $\kappa$.  In this case, the ordered sequence constructed in the proof of Lemma \ref{usefullemma} is actually of size $< \kappa'$.  This situation allows us to repeat the above proof and construct $2^\kappa$ non-isomorphic models of size $\kappa$.  Many other cardinal arithmetic set-ups suffice for many models.
		
\underline{$(E)$}\\
We have already mentioned that Existence for $\kappa$-coheir is equivalent to the fact every model is $\kappa$-Galois saturated.  The following theorem shows that this follows from categoricity in a $\kappa$-closed cardinal.
\begin{theorem} \label{catexist}  Suppose $K$ is an AEC satisfying the amalgamation property, JEP and has no maximal models.
If $K$ is categorical in a cardinal $\lambda$ satisfying $\lambda = \lambda^{< \kappa}$, then every member of $K_{\geq \chi}$ is $\kappa$-Galois saturated, where $\chi = min \{ \lambda, sup_{\mu < \kappa}( \beth_{(2^{\mu})^+})\}$.
\end{theorem}

{\bf Proof:}  First, note that by the assumptions on $K$ and   the assumption that $\lambda = \lambda^{< \kappa}$  we can construct a $\kappa$-Galois saturated member of $K_\lambda$.  Since this class is categorical, all members of $K_\lambda$ are $\kappa$-Galois saturated.
	
	The easy case is when $\lambda<\chi$:  Suppose $M \in K$ is not $\kappa$-saturated and $\|M\| > \lambda$.  Then there is some small $M^- \prec M$ and $p \in gS(M^-)$ such that $p$ is not realized in $M$.  Then let $N \prec M$ be any substructure of size $\lambda$ containing $M^-$.  Then $N$ doesn't realize $p$, which contradicts its $\kappa$-Galois saturation.
	
	For the hard part, suppose $M \in K$ is not $\kappa$-Galois saturated and $\|M\| \geq \sup_{\mu<\kappa}(\beth_{(2^{\mu})^+})$.  
	There is some small $M^- \prec M$ and $p \in gS(M^-)$ such that $p$ is not realized in $M$.  We define a new class $(K^+, \prec^+)$ that depends on $K, p $ and $M^-$ as follows:
	
 $L(K^+): = L(K) \cup \{ c_{m} : m \in |M^-|\}$ by
	\begin{eqnarray*}
	K^+ &=& \{ N : N \te{ is an }L(K^+)\te{ structure st } N\rest L(K) \in K,  \te{  there exists} \\
	&{}& h:M^-\rightarrow N\rest L(K) \te{ a }K\te{-embedding such that }
	h(m)=(c_{m})^N \\
	&{}& \te{ for all  }m\in M^- \te{ and }
 N\rest L(K) \te{ omits }h(p) \}.\\
	N_1 \prec^+ N_2 &\iff& N_1 \rest L(K) \prec N_2\rest L(K) \te{ and } N_1 \subset_{L(K^+)} N_2.
	\end{eqnarray*}
	This is clearly an AEC with $LS(K^+) = \|M^-\|< \kappa$ and $\seq{M, m}_{m \in |M^-|} \in K^+$.  
	
	By Shelah's Presentation Theorem $K^+$ is a $PC_{\mu,2^\mu}$ for $\mu:=LS(K^+)$.  By \cite[Theorem VII.5.5]{shelahfobook} the Hanf number of $K^+$ is $\leq\beth_{(2^\mu)^+}\leq\chi$.
	
	Thus, $K^+$ has arbitrarily large models.  In particular,  there exists $N^+\in K_\lambda^+$.  Then $N^+ \rest L(K) \in K_\lambda$ is not $\kappa$-Galois saturated as it omits its copy of $p$. \hfill \dag\\

Regarding Extension, the strength of this assumption is not entirely clear.  The first order version is proved with compactness and the fact that it holds under strongly compact cardinals (see Theorem \ref{sceasy}) does not work to separate from compactness.  However, it does indeed hold in nonelementary classes; see the discussion of quasiminimal classes and $\lambda$-saturated models below.  Very recently, Vasey \cite{vaseyc} is able to show that, if $\kappa = \beth_\kappa$ in addition to our other hypotheses and $K$ is $\lambda$-Galois stable for $\lambda \geq \kappa$ instead of no order property and $\chi$-tame for some $\chi < \kappa$, then Extension holds (and moreover that $\kappa$-coheir is an independence relation) for the class of $\lambda^+$-saturated models.  

\begin{remark}
While the rest of the results use that $K$ satisfies all of Hypothesis \ref{hyp}, the proof of Theorem \ref{catexist} only uses the amalgamation property and also avoid any use of tameness or type shortness.  Note that Propositions \ref{tameind-prop} and \ref{orderprop} and Theorem \ref{catexist} are not used in the rest of the paper, but are intended to motivate the hypotheses of Theorem \ref{mainthm} as natural.
\end{remark}

At the referee's request, we explain the difficulties in adapting the first-order proof to the current context.  The main difficulty is the standard one in trying to transfer results from first-order to AECs: the lack of syntax and compactness.  The lack of syntax was helped by the recent isolation of the properties of tameness and type shortness.  These properties allow us to treat Galois types as maximal collections of small Galois types, in the same way they are maximal collections of formulas in first-order.  However, this does not solve all problems as there are still two key differences: types are only well-behaved over models and types cannot be constructed inductively in an effective manner (i. e., every type is complete over it's domain).  Nonetheless, this intuition allows the authors to adapt many arguments from first-order.  Regarding compactness, tameness and type shortness are approximations, but assuming Extension outright is the strongest approximation.

Before continuing, we also identify a few contexts which are known to satisfy this hypothesis, especially (1), (2), and (3) of Theorem \ref{mainthm}.
\begin{itemize}
	\item {\bf First order theories:}  Since types are syntactic and over sets, they are $< \aleph_0$ tame and $< \aleph_0$ type short and (4) follows by compactness.  Additionally, when (3) holds, the theory is stable so coheirs are equivalent to non-forking.  %While we don't claim to have discovered anything new about first-order theories, formally speaking our framework apply to $K_T$ where $T$ is a superstable first-order theory and $K_T$ is the class of $|T|^+$-saturated models (our $\kappa$ is $|T|^+$).
	
	\item {\bf Large cardinals:}  Boney \cite{tamelc} proves that (1),(2), and Extension hold for any AEC $K$ that are essentially below a strongly compact cardinal $\kappa$ (this holds, for instance, if $LS(K) < \kappa$).  Slightly weaker (but still useful) versions of (1) and (2) also hold if $\kappa$ is measurable or weakly compact.  See Section \ref{lcrevisitedsection} for more.
	
	\item {\bf Homogeneous model theory:}  The sequential homogeneity of the monster model means that Galois types are syntactic, so we have $<\aleph_0$ tameness and type shortness as above.  See, for instance, 
%	Hyttinen and Shelah \cite{strsplithom} 
Grossberg and Lessmann \cite{GrLe1} 
	for more discussion and references.
	
	\item {\bf Quasiminimal classes:}  The quasiminimal closure operation means that Galois types are quantifier free types and amalgamation and other properties are proved in the course of proving categoricity; see \cite[p. 190]{baldwinbook}.

	\item {\bf Saturated models of a superstable theory:} Let $T$ be a superstable first-order theory and $K^\lambda_T$ be the class of $\lambda$-saturated models of $T$ ordered by elementary substructure.  This is a nonelementary class, but still satisfies our hypotheses.

	\item {\bf Averageable classes:} Averageable classes are $EC(T, \Gamma)$ classes that have a suitable ultraproduct-like relation (that averages the structures), see Boney \cite{gamult}.  Examples of averageable classes are dense ordered group with a cofinal $\mathbb{Z}$-chain and ordered vector fields.  The existence of an ultraproduct-like operation mean that $(E)$ (and much more) can be proved similar to the first order version.  Amongst these classes, the question of stability in the guise of no order property becomes crucial.  Torsion modules over a PID are an example that lie on the stable side of the line.

\end{itemize}

Since the circulation of early drafts, the notion of $\kappa$-coheir has been used and extended by various authors, especially Boney, Grossberg, Kolesnikov, and Vasey in \cite{bgkv} and Vasey in \cite{vaseybehemoth}, \cite{vaseyc} and \cite{gvas}.  Moreover, \cite[Theorem 6.7]{bgkv} has shown that, if $\kappa$-coheir is an independence relation, then it is the \emph{only} independence relation.

\begin{theorem}[\cite{bgkv}]
Under the hypotheses of Theorem \ref{mainthm}, $\dnf$ is the only independence relation on $K_{\geq \kappa}$.  In particular, if $\adnf$ satisfies $(I)$, $(M)$, $(E)$, and $(U)$, then $\adnf = \dnf$ on $K_{\geq \kappa}$.
\end{theorem}

%%%%%%%%%%%%%%%%%%%%%%%%%%%%%%%%%%%%%%%%%%%%%%%
\section{Getting Local Character} \label{lcsection}
%%%%%%%%%%%%%%%%%%%%%%%%%%%%%%%%%%%%%%%%%%%%%%%

Local character is a very important property for identifying dividing lines.  In the first order context, some of the main classes of theories--superstable, strictly stable, strictly simple, and unsimple--can be identified by the value of $\kappa(T)$.  By finding values for $\kappa_\alpha(\dnf)$ under different hypotheses, we get candidates for dividing lines in AECs.

Readers familiar with first order stability theory will recall that there is a notion of an heir of a type that is the dual notion to coheir, which our nonforking is based on.  Heir is equivalent to the notion of coheir under the assumption of no order; see \cite{pillay} as a reference.  We develop a Galois version of heir and show it is equivalent to heir under the assumption of no weak $\kappa$-order property.   This equivalence allows us to adapt an argument of \cite{shvi635} to calculate $\kappa^*_\omega(\dnf)$ from categoricity.  

Recall that we are working under Hypothesis \ref{hyp}.  This means that, unless stated in the hypothesis, we only have the properties of nonforking that follow immediately from the definition, like those in Theorem \ref{easypropthm}.  We explicitly state any other assumptions.  In particular, note that Theorem \ref{ulcthm} doesn't assume $(E)$, tameness, or type shortness.

Recall that `small' refers to objects of size $< \kappa$.

\begin{defin}\label{heir-def}
We say that $p \in gS^I(N)$ is an \emph{heir} over $M \prec N$ iff for all small $I_0 \subset I$, $M^- \prec M$, and $M^- \prec N^- \prec N$ (with $M^-$ possibly being empty), there is some $f:N^- \to_{M^-} M$ such that $f(p^{I_0} \rest N^-) \leq p$; that is, $f(p)^{I_0} \rest f(N^-) = p^{I_0} \rest f(N^-)$.  We also refer to this by saying $p$ is a heir of $p \rest M$.
\[
\xymatrix{
M \ar[r] & N\\
M^- \ar[r]\ar[u] & N^- \ar[ul]_f \ar[u]}
\]
\end{defin}

\begin{remark}
Note that the diagram above and the related diagrams in the proof of Theorem \ref{ulcthm} are not commutative diagrams; this would require that $f$ fix all of $N^-$, trivializing it.  Instead, this diagram is included to help visualize the definition.
\end{remark}

At first glance, this property seems very different from the first order version of heir.  However, if we follow the remark after Theorem \ref{fulltheorem}, we can think of restrictions of $p$ as formulas and small models as parameters.  Then, $M^-$ is a parameter from $M$, $N^-$ is a parameter from $N$, $f(N^-)$ is the parameter from $M$ that corresponds to $N^-$ (notice that it fixes $M^-$), and $f(p\rest N^-) \leq p$ witnesses that it the original formula $p \rest N^-$ is still in $p$ with a parameter from $M$.  Vasey has made this intuition explicit in \cite[Proposition 5.3]{vaseybehemoth} with the notion of Galois Morleyization.

If we restrict ourselves to models, then the notions of heir over and nonforking (coheir over) are dual with no additional assumptions.
\begin{prop} \label{heircoduality}
Suppose $M_0\prec M, N$.  Then $gtp(M/N)$ does not fork over $M_0$ iff $gtp(N/M)$is an heir over $M_0$.
\end{prop}

{\bf Proof:} First, suppose that $M \dnf_{M_0} N$ and let $a \in {}^{<\kappa}N$.  To show $gtp(N/M)$ is an heir over $M_0$, let $M_0^- \prec M_0$ and $M^- \prec M$ both be of size $< \kappa$ such that $M_0^- \prec M^-$.  Find $N^- \prec N$ of size $< \kappa$ containing $M_0^-$ and $a$.  By the definition of nonforking, $gtp(M^-/N^-)$ is realized in $M_0$.  This means that there is $g \in Aut_{N^-} \sea$ such that $g(M^-) \prec M_0$.  Set $f = g \rest M^-$.  Then $f: M^- \to_{M_0^-} M_0$ such that $f(gtp(a/M^-)) = gtp(a/f(M^-))$.  Since $a$, $M_0^-$, and $M^-$ were arbitrary, $gtp(N/M)$ is an heir over $M_0$, where $gtp(N/M)$, $N$, and $M_0$ here stand for $p$, $N$, and $M$ in Definition \ref{heir-def}.

Second, suppose that $gtp(N/M)$ is an heir over $M_0$.  Let $b \in M$ and $N^- \prec N$ both be of size $< \kappa$.  Since $M$ is a model, we may expand $b$ to a model $M^- \prec M$ of size $< \kappa$.  Then, if we can realize $gtp(M^-/N^-)$ in $M_0$, we can find a realization of $tp(b/N^-)$ there as well.  By assumption, there is some $f: M^- \to M_0$ such that $gtp(f(N^-)/f(M^-)) = gtp(N^-/f(M^-))$.  This type equality means that there is some $g \in Aut_{f(M^-)} \sea$ such that $g \circ f$ is the identity on $N^-$.  Thus, $g \circ f$ is in $Aut_{N^-} \sea$ and sends $M^-$ to $f(M^-) \prec M_0$.  Thus, $gtp(M^-/N^-) = gtp(g \circ f(M^-)/N^-)$ and is realized in $M_0$, as desired.  Since $b$ and $N^-$ were arbitrary, $gtp(M/N)$ does not fork over $M_0$. \hfill \dag\\

Proposition \ref{heircoduality} was proven just from the definitions, without assuming any tameness or type shortness.  If we assume even the weak symmetry $(S^*)$, then we have that nonforking and heiring are equivalent for models.  However, the goal is to show that they are equivalent for all Galois types.  Assuming full symmetry $(S)$ is enough to get this implication in one direction.

\begin{theorem}
Suppose $\dnf$ satisfies $(S)$.  If $p \in gS(N)$ and $M \prec N$, then $p$ does not fork over $M$ holds implies $p$ is an heir over $M$.
\end{theorem}

{\bf Proof:} Suppose $p \in gS(N)$ does not fork over $M$.  Then, given $A$ that realizes $p$, we have $A \dnf_M N$.  By $(S)$, we can find $M^+ \succ M$ containing $A$ such $N \dnf_M M^+$.  By Proposition \ref{heircoduality}, we then have $gtp(M^+/N)$ is an heir over $M$.  By monotonicity, $p = gtp(A/N)$ is an heir over $M$.\hfill \dag\\

For the other direction, we use the weak $\kappa$-order property (recall Definition \ref{orderdef}).

\begin{theorem} \label{order2heir}
Suppose $M \prec N$ and $M$ is $\kappa$-Galois saturated.  If there is $p \in gS(N)$ that is an heir over $M$ and not a coheir over $M$, then $K$ has the weak $\kappa$-order property and it is witnessed in $M$.
\end{theorem}

{\bf Proof:} It suffices to show this under the assumption that $|\ell(p)| < \kappa$.  Suppose $p = gtp(b/N)$ is a heir over $M$ but not a coheir.  Thus, there is some $N^- \prec N$ such that $gtp(b/N^-)$ is not realized in $M$.  Possibly by expanding $N^-$, we may assume that there is some small $M_0 \prec M$ such that $M_0 \prec N^-$\footnote{The introduction of $M_0$ is not necessary, but avoids added complication that would be necessary to introduce Galois types over the empty set}.  Because $p$ is a heir, there is some $f: N^- \to_{M_0} M$ such that $f( p \rest N^-) = p \rest f(N^-)$.  By $\kappa$-Galois saturation, there is some $b' \in M$ such that 
$$b' f(N^-) \vDash gtp(b f(N^-)/M_0) = gtp(f(b) f(N^-)/M_0) = gtp(b N^-/M_0)$$
Set $p_0 := gtp(N^- b/M_0)$ and $q_0 := gtp(N^- b'/M_0)$\footnote{In stating this, we have given the sets $N^-b$ and $N^-b'$ a ``natural'' enumeration: $N$ is listed first and in the same order in each enumeration; then, $b$ is listed in the order inherited from $p$ and $b'$ is inherited in the order inherited from realizing the same type as $b$ over $M_0$.  Note we have switched the ordering of $b$ and $N^-$ to make our conclusion match Definition \ref{orderdef}.}.  Note that these are different; we know that $N^- b'$ does not realize $p_0$ by assumption.  Now we will construct, by induction, sequences $a_i, b_i \in M$ such that
\begin{enumerate}
	\item $a_i b \vDash p_0$ for $i < \kappa$ (note that this implies $a_i \cong N^-$);
	\item $N^- b_j \vDash q_0$ for $i < \kappa$; and
	\item \label{second} \
	\begin{eqnarray*}
i \leq j \implies gtp(a_i b_j/M_0) = p_0\\
i > j \implies gtp(a_i b_j/M_0)= q_0
\end{eqnarray*}
\end{enumerate}
Note that (\ref{second}) immediately implies the weak $\kappa$-order property.

For $i = 0$, the elements are already built: $b_0 = b'$ and $a_0 = f(N^-)$.

Suppose $i < \kappa$ and $\{ a_j, b_j : j < i\}$ are already defined.  By induction, the following is true
\begin{enumerate}
\item $a_j b \vDash p_0$ for all $j < i$;
\item $N^- b_j \vDash q_0$ for all $j < i$; and
\item $N^- b \vDash p_0$.
\end{enumerate}
Take some $M^- \prec M$ that contains $\{ a_j, b_j : j < i\} \cup M_0$ and let $N^+ \prec N$ contain $M^-$ and $N^-$.  Because $p \in gS(N)$ is an heir over $M$, there is $g: N^+\to_{M^-} M$ such that $gtp(b/g(N^+)) = gtp(g(b)/g(N^+))$.  For $j < i$, $g$ fixes $M_0$ and $b_j$.  Thus,
$$gtp\left( g(N^-) b_j/M_0 \right) = gtp\left( g(N^-) g(b_j) / M_0 \right) = gtp(N^- b_j/M_0)$$
Additionally, we have that $gtp(b g(N^+)/M_0) = gtp(b N^+/M_0)$ and, thus, $gtp(a_i g(N^-)/M_0) = gtp(a_iN^-/M_0)$.  Setting $a_i := g(N^-)$, we have shown
\begin{enumerate}
	\item $a_j b \vDash p_0$ for all $j < i$;
	\item $a_i b_j \vDash q_0$ for all $j < i$; and
	\item $a_i b \vDash p_0$.
\end{enumerate}
Now take small $M^+ \prec M$ to contain $M^-$ and $a_i$.  By $\kappa$-Galois saturation, there is some $b_i \in M$ such that $gtp(b_i/M^+) = gtp(b/M^+)$.  Because $M^+$ contains $\{a_j : j \leq i\}$, this means that
\begin{enumerate}
	\item $a_j b_i \vDash p_0$ for all $j < i$;
	\item $a_i b_j \vDash q_0$ for all $j < i$; and 
	\item $a_i b_i \vDash p_0$.
\end{enumerate}
This completes the induction and the proof.\hfill \dag\\

Now that we have established an equivalence between nonforking and being an heir, we aim to derive local character.  For this, we use heavily the proof of \cite[Theorem 2.2.1]{shvi635}, which shows that, under certain assumptions, the universal local character cardinal for non-splitting is $\omega$ (recall Definition \ref{lc-def}).  Examining the proof,  much of the work is done by basic independence properties--namely (I), (M), and (T)--and the other assumptions on $K$--namely categoricity, amalgamation, and EM models, which follow from no maximal models.  Only in case (c), defined below, do they need the exact definition of their independence relation (non $\mu$-splitting) and GCH.  In this case, we can use the definition of heir to complete the proof without any appeal to cardinal arithmetic.

At the referee's insistence, we include a complete proof.  This gives us an opportunity to expand on the proof given in \cite{shvi635} and correct some mistakes (e. g., they cite a different club guessing principle than the one they use).  We would like to thank Monica VanDieren and Sebastien Vasey in helping bridge some of the gaps.  A write up of a detailed proof for common generalization of \cite{shvi635} and Theorem \ref{ulcthm} can be found in \cite{svthm}.  Note that the argument below assumes $\kappa$ regular, but \cite{svthm} shows how to remove that assumption.

The proof uses the notion of limit models.  
\begin{defin}
Let $M \in K_\lambda$ and $\alpha < \lambda^+$ be a limit ordinal.  $N$ is $(\lambda, \alpha)$-limit over $M$ iff there is a resolution $\seq{N_i\in K_\lambda : i < \alpha}$ of $N$  so $N_0 = M$ and $M_{i+1}$ is universal over $M_i$.
\end{defin}
Limit models have been suggested as a substitute for Galois saturated models and the question of uniqueness of limit models has been suggested as a dividing line for AECs; see Shelah \cite{sh576}.  We discuss the uniqueness of limit models after the proof, but the key use of uniqueness of limit models in this proof are the following basic facts.

\begin{fact}\label{lm-fact}\
\begin{enumerate}
	\item If $M^\ell$ is $(\mu, \alpha_\ell)$-limit and $\cf \alpha_1 = \cf \alpha_2$, then $M^1 \cong M^2$.  Moreover, if both are limit over $M_0$, then $M^1 \cong_{M_0} M^2$.
	\item If $\seq{M^\ell_j \mid j \leq \delta}$ are increasing continuous chains for $\ell=1,2$ and $M^\ell_{j+1}$ is $(\lambda, \beta)$-limit over $M^\ell_j$ for all $j < \delta$, then there is $f:M^1_\delta \cong M^2_\delta$ such that $f(M^1_j) = M^2_j$ for each $j \leq \delta$.
\end{enumerate}
\end{fact}

{\bf Proof sketch:} The first is proved by a back and forth and the second is a consequence of the first. \hfill \dag\\

\begin{theorem}\label{ulcthm}
Let $\kappa \leq \mu < \lambda$ and suppose $\kappa$ is regular\footnote{Note that we don't make this assumption elsewhere.%\\WB: I'm not convinced this is necessary.  At the very least, we should be able to remove regularity at the cost of requiring that $\mu > \kappa$ and replacing $(\mu, \kappa)$-limit (which give us our $\kappa$-saturation) with $(\mu, \kappa^+)$-limit
}.  Suppose that $K$ has no weak $\kappa$-order property and is categorical in $\lambda$.  Then $\kappa_\omega^*(\dnf) = \omega$.  That is, for each $\alpha = \cf \alpha < \mu^+$, if
\begin{enumerate}

	\item $\seq{M_i \in K_\mu : i \leq \alpha}$ is increasing and continuous;
	
	\item each $M_{i+1}$ is universal over $M_i$; and
		
	\item $p \in gS^{< \omega}(M_\alpha)$

\end{enumerate}
then, for some $i < \alpha$, $p$ does not fork over $M_i$
\end{theorem}

{\bf Proof:} As mentioned above, we follow the proof of \cite[Theorem 2.2.1]{shvi635}.  By \cite[Theorem 8.21]{baldwinbook}, we have Galois stability in $\mu$ and $\kappa$.  

We want to show that any failure of the theorem implies there is a failure that takes one of three possible forms:
\begin{enumerate}

	\item[(a)] There exists a sequence $\seq{M_i : i \leq \alpha}$ and $p$ as in the hypothesis with $M_{i+1}$ being $(\mu, \kappa)$-limit over $M_i$ such that for all $i < \alpha$, $p \rest M_i$ does not fork over $M_0$ and $p$ forks over $M_i$.%\footnote{WB: It feels like the proof of this case occurs as a subcase of part (b) and we really worry about the case where $p \rest M_{2i+2}$ does not fork over $M_{2i+1}$.}
	\item[(b)] There exists a sequence $\seq{M_i : i \leq \alpha}$ and $p$ as in the hypothesis with $M_{i+1}$ being $(\mu, \kappa)$-limit over $M_i$ and $\kappa \leq \alpha$ such that for all $i < \alpha$, $p \rest M_{i + 1}$ and $p$ fork over $M_i$.
	\item[(c)] There exists a sequence $\seq{M_i : i \leq \alpha}$ and $p$ as in the hypothesis with $M_{i+1}$ being $(\mu, \kappa)$-limit over $M_i$ and $\alpha < \kappa$ such that, for all $i < \alpha$, $p \rest M_{2i+1}$ forks over $M_{2i}$, $p \rest M_{2i+2}$ does not fork over $M_{2i+1}$, and $p$ forks over $M_i$.

\end{enumerate}

Note that $(\mu, \kappa)$-limit implies $\kappa$-saturation and the crucial restriction on $\alpha$ in (b) and (c).  Additionally, in cases (b) and (c), the condition that $p$ forks over each $M_i$ follows from other forking in the case, so doesn't need to be guaranteed during the construction.

\begin{claim} \label{firstclaim}
If (a) and (b) are false, then for every $\seq{M_i : i \leq \alpha}$ and $p$ as in the hypothesis with $\alpha \geq \kappa$, there is $i < \alpha$ such that $p$ does not fork over $M_i$.
\end{claim}

The assumption that $\seq{M_i : i \leq \alpha}$ and $p$ are as in the hypothesis includes that $\alpha$ is regular.\\

{\bf Proof:}  We do a proof by contradiction.  Let $\seq{M_i \mid i \leq \alpha}$ and $p$ be as in the hypothesis with $p$ forking over each $M_i$.  First, using amalgamation and stability in $\mu$, we can assume that each $M_{i+1}$ is $(\mu, \kappa)$-limit over $M_i$.\\

{\bf Subclaim:} If (a) is false, then for any $\seq{N_i \mid i \leq \alpha}$ and $q$ as in the hyptohesis with $q$ forking over each $N_i$, there is an increasing, continuous sequence $\seq{j_i \mid i \leq \alpha}$ of ordinals such that $q \rest N_{j_{i+1}}$ forks over $N_{j_i}$.\\

{\bf Proof:} For each $i<\alpha$, consider the sequence $\seq{N_j \mid i \leq j \leq \alpha}$ starting at $i$. Since (a) fails, there is some $j(i) >i$ such that $q \rest N_{j(i)}$ forks over $N_i$. Form $j_i$ by setting $j_0 = 0$, taking unions at limits, and setting $j_{i+1} = j(j_i)$.  Then $q \rest N_{j_{i+1}}$ forks over $N_{j_i}$, as desired.\hfill$\dag_{\text{Subclaim}}$\\

This gives us an increasing sequence $\seq{j_i \mid i \leq \alpha}$.  Define an increasing, continuous $\seq{M^*_i \mid i \leq \alpha}$ by $M^*_i = M_{j_i}$.  Since amalgamation holds, $M^*_{i+1}$ is still $(\mu, \kappa)$-limit over $M_i^*$.  Thus, $\seq{M^*_i\mid i \leq \alpha}$ is a sequence that the hypothesis says does not exist. \hfill $\dag_{\text{Claim \ref{firstclaim}}}$\\

%The hypotheses of Claim \ref{firstclaim} hold for any regular $\alpha \geq \kappa$ when (a) and (b) are false.  This immediately gives that, if (a) and (b) are false, then the theorem holds for sequences of length at least $\kappa$ (that is, $\kappa^*_\omega(\dnf) \leq \kappa$).

\begin{claim} \label{secondclaim}
If the theorem is false, then at least one of (a), (b), and (c) holds.
\end{claim}

{\bf Proof:}  Suppose that $\seq{M_i \mid i \leq \alpha}$ and $p$ are in the hypothesis such that, for all $i < \alpha$, $p$ forks over $M_i$ and that (a) and (b) are false.  Claim \ref{firstclaim} shows that the failure of (a) and (b) means any such counterexample must have $\alpha < \kappa$.  Again, with amalgamation and stability, we can assume that $M_{i+1}$ is $(\mu, \kappa)$-limit over $M_i$.  Let $\seq{M_{i, j} \mid j \leq \kappa}$ witness this.  Since (a) is false, by the proof of Claim \ref{firstclaim}, without loss of generality, $p \rest M_{i+1}$ forks over $M_i$.  Since we know the theorem holds for $\kappa$ length sequences, there is some $j_i < \kappa$ such that $p \rest M_{i+1}$ does not fork over $M_{i, j_i}$. %By (M), we may assume that $\cf j_i = \kappa$ by replacing $j_i$ with $j_i + \kappa$.  
Also by (M), we have that $p \rest M_{i+1, j_i}$ forks over $M_i$.  Then define increasing continuous $\seq{M^*_i \mid i \leq \alpha}$ by
\begin{eqnarray*}
M^*_i = \begin{cases}
M_{2k} & \text{ if } i = 2k\\
M_{2k+1, j_{2k+1}} & \text{ if } i = 2k+1
\end{cases}
\end{eqnarray*}
This gives a sequence as in (c).\hfill $\dag_{\text{Claim \ref{secondclaim}}}$\\

Now we show that each of (a), (b), and (c) is false, which will finish the proof.  We examine case (b) first, as it is the most different from the argument in \cite{shvi635}.  Then, we examine cases (a) and (c) together as the initial construction is the same in each.\\

{\bf Case (b):} Suppose that $\seq{M_i : i \leq \alpha}$ and $p$ are as in (b).  Find the minimum $\sigma$ such that $2^\sigma > \kappa$.  Then $\sigma \leq \kappa$ and $2^{< \sigma} \leq \kappa$.  By using the contrapositive of (M), we have that $p \rest M_\sigma$ forks over $M_i$ for all $i < \sigma$.  Thus, $\seq{M_i : i \leq \sigma}$ and $p \rest M_\sigma$ satisfy the forking relations of (b), although we don't know that $\sigma$ is regular.

Set $M = M_\sigma$.  By Theorem \ref{order2heir} and the assumption of no weak $\kappa$-order property, we know that $p \rest M_{i+1}$ is not a heir over $M_i$ for all $i < \sigma$.  We are going to contradict stability in $\kappa$ by finding $2^\sigma$ many types over a model of size $2^{< \sigma}$.\\

{\bf Construction 1:} There are $\seq{M^i, N^i \in K_{<\kappa} \mid i < \sigma}$ such that
\begin{enumerate}
	\item $M^i \prec M_i$ and $N^i \prec M_{i+1}$; and
	\item if $f:N^i \to_{M^i} M_i$, then $f(p \rest N^i) \neq p \rest f(N^i)$.
\end{enumerate}
This follows directly from the hypothesis: for each $i < \sigma$, $p \rest M_{i+1}$ is not a heir over $M_i$.  Negating Definition \ref{heir-def} gives us the models.\\

{\bf Construction 2:} There are models $\{\hat M_i, \hat N_i \in K_{<\kappa} \mid i < \sigma\}$\footnote{Note we have truncated the chain to length $\sigma$.  We could continue this construction all the way to $\kappa$, but going beyond this would require greater saturation.} and $K$-embeddings $\{g_i \mid i < \sigma\}$ such that
\begin{enumerate}
	\item $\seq{\hat M_i \mid i < \sigma}$ is increasing;
	\item $M^{i+1} \prec \hat M_i \prec M_{i+1}$ and $\hat M_i \prec \hat N_i \prec \hat M_{i+1}$;
	\item $N^{i+1} \prec \hat N_i$; and
	\item $g_i:\hat N_i \to_{\hat M_i} M_{i+1}$ with $g_i(\hat N_i) \prec \hat M_{i+1}$.
\end{enumerate}

We build this by induction.  Set $\hat M_0 := M_0$.  If $i$ is limit, let $\hat M_i \prec M_{i+1}$ be of size $< \kappa$ that contains $M^{i+1}$ and $\cup_{j < i} \hat M_j$.

Given $\hat M_i$, we wish to build $\hat N_i$, $g_i$, and $\hat M_{i+1}$ (this will complete the construction).  Let $\hat N_i \prec M_{i+2}$ be of size $<\kappa$ that contains $\hat M_i$ and $N^{i+1}$.  $M_{i+1}$ is $\kappa$-Galois saturated, so it realizes $gtp(\hat N_i/\hat M_i)$.  This is witnessed by some $g_i: \hat N_i \to M_{i+1}$.  Now let $\hat M_{i+1} \prec M_{i+2}$ be of size $<\kappa$ that contains $\hat N_i$ and $g_i(\hat N_i)$.  This completes the construction.\\

{\bf Construction 3:} There are increasing $\seq{N_\eta \in K_{<\kappa} \mid \eta \in {}^{<\sigma} 2}$ and increasing isomorphisms $\{h_\eta \mid \eta \in {}^{<\sigma} 2\}$ such that
\begin{enumerate}
	\item $h_\eta: \hat M_{\ell(\eta)} \cong N_\eta$; and 
	\item $h_{\eta^\frown\seq{0}} \circ g_{\ell(\eta)}(\hat N_{\ell(\eta)}) \prec N_{\eta ^\frown \seq{0}} \cap N_{\eta^\frown \seq{1}}$.
\end{enumerate}

We build this by induction on the length of $\eta$. For the base case $\eta = \emptyset$, set $N_\eta = \hat M_0$ and $h_\eta$ to be the identity.  If $\ell(\eta)$ is limit, then we have $h_\eta^- = \cup_{\alpha < \ell(\eta)} h_{\eta \rest \alpha}$ is an isomorphism from $\cup_{\alpha < \ell(\eta)} \hat M_\alpha$ to $\cup_{\alpha < \ell(\eta)} N_{\eta \rest \alpha}$.  We have $\cup_{\alpha < \ell(\eta)} \hat M_\alpha \prec M_{\ell(\eta)}$, so we can extend $h_\eta^-$ to an isomorphism with domain $M_{\ell(\eta)}$.  Set $N_\eta := h_\eta(M_{\ell(\eta)})$.

Suppose we are given $N_\eta$ and $h_\eta$. Since we are working inside the monster model $\sea$, we can extend $h_\eta$ to $h_\eta^+ \in Aut \sea$ and $g_{\ell(\eta)}$ to $g_{\ell(\eta)}^+ \in Aut \sea$.  Then set
\begin{itemize}
	\item $h_{\eta^\frown\seq{0}} := h_\eta^+ \rest \hat M_{\ell(\eta)+1}$;
	\item $h_{\eta^\frown\seq{1}} := h_\eta^+ \circ g^+_{\ell(\eta)} \rest \hat M_{\ell(\eta)+1}$; and 
	\item $N_{\eta^\frown\seq{i}} := h_{\eta\frown\seq{i}} ( M_{\ell(\eta)+1})$ for $i = 0,1$.
\end{itemize}
Since $g_{\ell(\eta)}$ fixes $\hat M_{\ell(\eta)}$, $h_{\eta^\frown\seq{1}}$ extends $h_\eta$.
Now, we have $\hat N_{\ell(\eta)} \prec \hat M_{\ell(\eta) + 1}$, so $h_{\eta^\frown\seq{0}} \circ g_{\ell(\eta)}(\hat N_{\ell(\eta)}) \prec N_{\eta^\frown \seq{1}}$ by applying $h_{\eta^\frown\seq{1}}$ to both sides.  Similarly, we have $g_{\ell(\eta)}(\hat N_{\ell(\eta)}) \prec \hat M_{\ell(\eta)+1}$, so $h_{\eta^\frown\seq{0}} \circ g_{\ell(\eta)}(\hat N_{\ell(\eta)}) \prec N_{\eta^\frown \seq{0}}$ by applying $h_{\eta^\frown\seq{0}}$ to both sides.\\

Now that we have completed our constructions, define the following: 
\begin{itemize}
	\item $\hat M_\sigma = \cup_{i < \sigma} \hat M_i$; 
	\item for $\eta \in {}^\sigma 2$, set $h_\eta := \cup_{\alpha < \sigma} h_{\eta \rest \alpha}$ and $N_\eta = \cup_{\alpha < \sigma} N_{\eta \rest \alpha}$ so $h_\eta: \hat M_\sigma \cong N_\eta$; and
	\item $N_* \prec \sea$ of size $\kappa$ to contain $\cup_{\eta \in {}^{<\sigma} } N_\eta$.
\end{itemize}
For each $\eta \in {}^\sigma 2$, we have $h_\eta(N_\eta) \prec N_*$.  Let $p_\eta$ be some extension of the type $h_\eta(p \rest \hat M_\sigma)$ to $N_*$.

\begin{claim}\label{manytypesclaim}
 If $\eta \neq \nu \in {}^\sigma 2$, then $p_\eta \neq p_\nu$.
\end{claim}

Set $\rho = \eta \cap \nu$ and WLOG assume $\rho^\frown\seq{0} \subset \eta$ and $\rho^\frown\seq{1} \subset \nu$.  Set $i = \ell(\rho)$.  We built $g_i:\hat N_i \to_{\hat M_i} M_{i+1}$ with $N^{i+1} \prec \hat N_i$ and $M^{i+1} \prec \hat M_i$ in Construction 2.  By Construction 1, this implies 
$$g_i(p \rest N^{i+1}) \neq p \rest g_i(N^{i+1})$$
Inequality of Galois types transfers up, so
$$g_i(p \rest \hat N_i) \neq p \rest g_i(\hat N_i)$$
Applying $h_{\rho^\frown\seq{0}}$ to both sides, we get
$$h_{\rho^\frown\seq{0}} \circ g_i( p \rest \hat N_i) \neq  h_{\rho^\frown\seq{0}}(p\rest  g_i(\hat N_i))$$

From Construction 3, we have that $h_{\rho^\frown\seq{1}} \rest \hat N_i = h_{\rho^\frown \seq{0}} \circ g_i \rest \hat N_i$.  Thus, we can rewrite the above as 
$$h_{\rho^\frown\seq{1}}( p \rest \hat N_i) \neq  h_{\rho^\frown\seq{0}}(p\rest  g_i(\hat N_i))$$
Then
\begin{itemize}
	\item $p \rest \hat N_i \leq p \rest \hat M_\sigma$ and $h_\nu \rest \hat N_i = h_{\rho^\frown\seq{1}} \rest \hat N_i$ by construction, so $h_{\rho^\frown \seq{1}}(p \rest \hat N_i) \leq p_\nu$; and, similarly,
	\item $p \rest g_i(\hat N_i) \leq p \rest \hat M_\sigma$ and $h_\eta \rest g_i(\hat N_i) = h_{\rho^\frown\seq{0}} \rest g_i(N_i)$, so $h_{\rho^\frown\seq{0}}(p \rest g_i(\hat N_i)) \leq p_\eta$.
\end{itemize}
Thus, $p_\nu \neq p_\eta$, as desired.\hfill $\dag_{\text{Claim \ref{manytypesclaim}}}$\\

Thus, we have constructed $2^\sigma > \kappa$ many Galois types over a model of size $\kappa$, contradicting Galois stability in $\kappa$. \hfill $\dag_{\text{Case (b)}}$\\

{\bf Cases (a) and (c):} Suppose we have $\seq{M_i \mid i \leq \alpha}$ and $p$ as in (a) or (c).  Recall the notation $S^{\mu^+}_\alpha = \{ \delta < \mu^+ \mid \cf \delta = \alpha\}$.  We say that $\bar C = \seq{C_\delta \mid \delta \in S^{\mu^+}_\alpha}$ is a $S^{\mu^+}_\alpha$-club sequence if each $C_\delta \subset \delta$ is club.  It is easy to see that club sequences exist (this will be enough for case (a)), and Shelah \cite{shg} proves the existence of club guessing club sequences in ZFC under various hypotheses.  We will describe a construction of a sequence of models $\bar N(\bar C)$ based on a club sequence and then plug in the necessary club sequence in each case.

Given a $S^{\mu^+}_\alpha$-club sequence $\bar C$, enumerate each club $C_\delta \cup \{\delta\}$ in increasing order as $\seq{\beta_{\delta, j} \mid j \leq \alpha}$; $C_\delta$ and the enumerations depend on the club sequence $\bar C$ and should be adorned with a $\bar C$ to indicate this dependence, but we will not do this.

\begin{claim} \label{lm-chain-claim}
We can build increasing, continuous $\bar N(\bar C) = \seq{N_i \in K_\mu \mid i < \mu^+}$ such that 
\begin{enumerate}
	\item $N_{i+1}$ is $(\mu, \kappa)$-limit over $N_i$;
	\item when $i \in S^{\mu^+}_\alpha$, there is $g_i:M_\alpha \cong N_i$ such that $g_i(M_j) = N_{\beta_{i, j}}$ for all $j \leq \alpha$; and
	\item when $i \in S^{\mu^+}_\alpha$, there is $a_i \in N_{i+1}$ that realizes $g_i(p)$.
\end{enumerate}
\end{claim}

The proof uses Fact \ref{lm-fact}.  Build an increasing, continuous chain of models such that successors are $(\mu, \kappa)$-limit over their predecessors and take limits at unions.  At successors of limits $i$ of cofinality $\alpha$, use Fact \ref{lm-fact} to build an isomorphism $g_i$ between $\seq{M_j \mid j \leq \alpha}$ and $\seq{N_{\beta_{i, j}} \mid j \leq \alpha}$ as described in (2).  For (3), find $N \in K_\mu$ and $b \in N$ such that $p = gtp(b/M_\alpha)$ and let $g_i^+ \in Aut \sea$ extend $g_i$.  Since $N_{i+1}$ is universal over $N_i$, there is some $f_i:g_i^+(N) \to_{N_i} N_{i+1}$.  Since $N_i = g_i(M)$, $a_i := f_i(g_i^+(b)) \in N_{i+1}$ realizes $g_i(p)$. \hfill $\dag_{\text{Claim \ref{lm-chain-claim}}}$\\

Because there are arbitrarily large models, by \cite[Theorem 8.18]{baldwinbook}, there is $\Phi$ that is proper for $K$ with $|\Phi| = LS(K) < \kappa$.  Our crucial use of categoricity (other than deriving stability) is that we may assume that $N:= \cup_{i < \mu^+} N_i \prec EM_\tau(\mu^+, \Phi)$.  Thus, we can write $a_i = \tau_i(\gamma_1^i, \dots, \gamma_{n(i)}^i)$ with $\gamma^i_1< \dots <\gamma^i_{m(i)} < i \leq \gamma^i_{m(i)+1} < \dots< \gamma^i_{n(i)} < \mu^+$.

  Now we begin to differentiate between the two cases.  In each, we will find $i_1 < i_2 \in S^{\mu^+}_\alpha$ such that $gtp(a_{i_1}/N_{i_1})$ and $gtp(a_{i_2}/N_{i_1})$ are both the same (because of the $EM$ structure) and different (because they exhibit different forking behavior), which is our contradiction.\\

{\bf Case (a):} Let $\bar{C}$ be an $S^{\mu^+}_\alpha$-club sequence, and set $\seq{N_i \in K_\mu \mid i < \mu^+} = \bar N(\bar C)$ from Claim \ref{lm-chain-claim}.  There is a stationary subset $S^* \subset S^{\mu^+}_\alpha$ such that
\begin{enumerate}
	\item for every $i \in S^*$, we have $\tau_i = \tau_*$; $n(i) = n_*$; $m(i) = m_*$; $\gamma^i_j = \gamma^*_j$ for $j \leq m_*$; and $\beta_{i,0} = \beta_{*,0}$.
\end{enumerate}

Set $E = \{\delta < \mu^+ \mid \delta$ is limit and $EM_\tau(\delta, \Phi) \cap N = N_\delta\}$.  This is a club.  Let $i_1 < i_2 \in S^* \cap E$.  Then we have
\begin{eqnarray*}
gtp\left(a_{i_1}/N_{i_1}\right) = gtp\left(\tau_*(\gamma^*_1, \dots,\gamma^*_{m_*}, \gamma^{i_1}_{m_*+1}, \dots, \gamma^{i_1}_{n_*})/N \cap EM_\tau(i_1, \Phi)\right)\\
=gtp\left(\tau_*(\gamma^*_1, \dots, \gamma^*_{m_*}, \gamma^{i_2}_{m_*+1}, \dots, \gamma^{i_2}_{n_*})/N \cap EM_\tau(i_1, \Phi)\right) = gtp\left(a_{i_2}/N_{i_1}\right)
\end{eqnarray*}
because the only difference in the two lie entirely above $i_1$.\\

We have that $g_{i_1}:(M_\alpha, M_0) \cong (N_{i_1}, N_{\beta_{*,0}})$ and that $p$ forks over $M_0$.  Thus, $gtp(a_{i_1}/N_{i_1}) = g_{i_1}(p)$ forks over $N_{\beta_{*, 0}}$.  On the other hand, $C_{i_2}$ is cofinal in $i_2$, so there is $j < \alpha$ such that $\beta_{i_2, j} > i_1$ and, thus, $N_{i_1} \prec N_{\beta_{i_2, j}}$.  Again, $g_{i_2}:(M_j, M_0) \cong (N_{\beta_{i_2, j}}, N_{\beta_{*,0}})$ and $p \rest M_j$ does not fork over $M_0$ because we are in case (a).  Thus, $gtp(a_{i_2}/N_{\beta_{i_2,j}}) = g_{i_2} (p \rest M_j)$ does not fork over $N_{\beta_{*, 0}}$.  By (M), $gtp(a_{i_2}/N_{i_1})$ does not fork over $N_{\beta_{*,0}}$.  Thus, $gtp(a_{i_1}/N_{i_1}) \neq gtp(a_{i_2}/N_{i_2})$, a contradiction.\hfill $\dag_{\text{Case (a)}}$\\

{\bf Case (c):}  Let $\chi$ be a big enough cardinal and create an increasing, continuous elementary chain $\seq{\mathfrak{B}_i \mid i < \mu^+}$ such that
\begin{enumerate}
	\item $\mathfrak{B}_i \prec (H(\chi), \in)$;
	\item $\|\mathfrak{B}_i\| = \mu$;
	\item \label{things} $\mathfrak{B}_0$ contains, as elements\footnote{When we say that $\mathfrak{B}_0$ contains a sequence as an element, we mean that it contains the function that maps an index to its sequence element.}, $\Phi$, $EM(\lambda^+,\Phi)$, $h$, $\mu^+$, $\seq{N_i \mid i < \mu^+}$, $S^{\mu^+}_\alpha$, $\seq{a_i \mid i \in S^{\mu^+}_\alpha}$,  and each $f \in \tau(\Phi)$; and %\footnote{WB: Add some comment about how we're coding everything low enough; maybe it doesn't matter}
%	\item $\seq{\mathfrak{B}_j \mid j \leq i} \in \mathfrak{B}_{i+1}$\footnote{WB: Pretty sure this is irrelevant};
	\item $\mathfrak{B}_i \cap \mu^+$ is an ordinal.
\end{enumerate}

In this case, we use the following:

\begin{fact}[\cite{shg}.Claim 2.3] \label{club-guess-fact-1}
Let $\lambda$ be a cardinal such that $\cf \lambda \geq \chi^{++}$ for some regular $\chi$ and let $S \subset S^\lambda_\chi$ be stationary.  Then there is a $S$-club sequence $\seq{C_\delta \mid \delta \in S}$ such that, if $E \subset \lambda$ is club, then there are stationarily many $\delta \in S$ such that $C_\delta \subset E$.
\end{fact}

%Below, a naked i that makes no sense should be a i_2 and a naked i' should be a i_1

We have that $\alpha < \kappa < \mu^+$, so there is an $S^{\mu^+}_\alpha$-club sequence $\bar C$ as in Fact \ref{club-guess-fact-1}.  Set $\seq{N_i \in K_\mu \mid i < \mu^+} = \bar N(\bar C)$ from Claim \ref{lm-chain-claim}.  Note that $E = \{ i < \mu^+ \mid \mathfrak{B}_i \cap \mu^+ = i\}$ is a club.  By Fact \ref{club-guess-fact-1}, there is some $i_2 \in S$ such that $C_{i_2} \subset E$.  We have 
$$a_{i_2} = \tau_{i_2}(\gamma_1^{i_2}, \dots, \gamma_{n(i_2)}^{i_2})$$
with $\gamma_1^{i_2} < \dots < \gamma_{m(i_2)}^{i_2} < i_2 \leq \gamma^{i_2}_{m(i_2)+1} < \dots < \gamma^{i_2}_{n(i_2)}$.  Since the $\beta_{i_2, j}$ enumerates a cofinal sequence in $i_2$, we can find $j < \alpha$ such that $\gamma_{m(i_2)}^{i_2} < \beta_{i_2, 2j+1} < i$.  Recall that we have $p \rest M_{2j+2}$ does not fork over $M_{2j+1}$ by the case hypothesis.  Then $(H(\chi), \in)$ satisfies the following formulas with parameters exactly the things listed in item (\ref{things}) above and ordinals below $\beta_{i_2, 2j+2}$
\begin{eqnarray*}
\exists x, y_{m(i_2)+1}, \dots, y_{n(i)}. ( ``x \in S" \wedge ``x > \beta_{i_2, 2j+1}" \wedge ``y_k \in (x, \mu^+) \text{ are increasing ordinals}"\\
\wedge ``a_x = \tau_{i_2}(\gamma_1^{i_2}, \dots, \gamma_{m(i_2)}^{i_2}, y_{m(i_2)+1}, \dots, y_{n(i_2)})" \wedge ``N_x \subset EM(x, \Phi)")
\end{eqnarray*}
witnessed by $x = i_2$ and $y_k = \gamma^{i_2}_k$.  By elementarity, $\mathfrak{B}_{\beta_{i_2, 2j+2}}$ thinks this as it contains all the parameters.  Let $i_1 \in (\beta_{i_2, 2j+1}, \mu^+) \cap \mathfrak{B}_{\beta_{i_2, 2j+2}} = (\beta_{i_2, 2j+1}, \beta_{i_2, 2j+2})$\footnote{This equality holds because $\beta_{i_2, 2j+2} \in C_{i_2} \subset E$ and is the key use of club guessing.} witness this, along with $\gamma'_{m(i_2)+1} < \dots < \gamma'_{n(i_2)} < \mu^+$.  Then we have
$$a_{i_1} = \tau_{i_2}(\gamma_1^{i_2}, \dots, \gamma_{m(i_2)}^{i_2}, \gamma'_{m(i_2)+1}, \dots, \gamma'_{n(i_2)})$$
with $\beta_{i_2, 2j+1} < \gamma_{m(i_2)+1}$.  We want to compare $gtp(a_{i_2}/N_{i_1})$ and $gtp(a_{i_1}/N_{i_1})$.
\begin{itemize}

	\item From the elementarity, we get that $N_{i_1} \subset EM(i_1, \Phi)$.  We also know that $i_1 < \beta_{i_2, 2j+2} < \gamma^{i_2}_{m(i_2)+1}, \gamma'_{m(i_2)+1}$.  Thus, as before, the types are equal.

	\item We know that $p \rest M_{2j+2}$ does not fork over $M_{2j+1}$.  Thus, $gtp(a_{i_2}/N_{\beta_{i_2, 2j+2}})$ does not fork over $N_{\beta_{i_2, 2j+1}}$.  Since we have $N_{\beta_{i_2, 2j+1}} \prec N_{i_1} \precneqq N_{\beta_{i_2, 2j+2}}$, this gives $gtp(a_{i_2}/N_{i_1})$ does not fork over $N_{\beta_{i_2, 2j+1}}$.
	
	\item We have $\beta_{i_2, 2j+1} < i_1$, so there is some $k < \alpha$ such that $\beta_{i_2, 2j+1} < \beta_{i_1, k} < i'$.  By assumption, $p$ forks over $M_k$.  Thus $g_{i_1}(p)$ forks over $N_{\beta_{i_1,k}}$.  Thus, $gtp(a_{i_1}/N_{i_1})$ forks over $N_{\beta_{i_2, 2j+1}} \prec N_{\beta_{i_1,k}}$.

\end{itemize}
As before, these three statements contradict each other.
\hfill $\dag_{\text{Case (c), Theorem \ref{ulcthm}}}$  \\

\begin{remark}
In the hypotheses of Theorem \ref{ulcthm}, we can replace the regularity of $\kappa$ with the assumption that $\mu > \kappa$.  To prove this statement, we replace $\kappa$ with $\kappa^+$ in several instances, e. g., making models $(\mu, \kappa^+)$-limit rather than $(\mu, \kappa)$-limit.
\end{remark}

%This construction could not go further than $\kappa$ many steps because the definition of heir requires all of the models and tuples involved to be of size $< \kappa$.  Thus, we need to know that stability fails at $\kappa$.  If we knew that nonforking and nonsplitting were the same, instead of just nonforking and heiring, then we would have a more general argument.  The connection between these two notions of independence and other is explored more in \cite{bgkv}.\footnote{WB: This no longer seems appropriate}

Once we have the universal local character, we can get results on the uniqueness of limit models.    Shelah and Villaveces \cite{shvi635} claimed uniqueness of limit models from categoricity.  However, VanDieren discovered a gap in the proof that the authors have not fixed.  VanDieren \cite{vandierennomax} \cite{nomaxerrata} proved uniqueness of limit models from categoricity with weaker assumptions than we have here, namely instead of full amalgamation it was assumed that only unions of limit models are amalgamation bases. Further results can be found in Grossberg, VanDieren, and Villaveces \cite{gvv}.

%Once we have the universal local character, we can get results on the uniqueness of limit models.  Limit models (introduced as brimmed in \cite{sh600}, see the definition below) have been suggested as a substitute for saturated models and the question of uniqueness has been suggested as a dividing line for AECs.  See Shelah \cite{sh576}, Shelah and Villaveces \cite{shvi635} claimed uniqueness of limit models from categoricity, however VanDieren in 1998 discovered a gap in the proof. In  2000 Shelah admitted that he is unable to fix that gap.   VanDieren \cite{vandierennomax} \cite{nomaxerrata} proved uniqueness of limit models from categoricity with weaker assumptions than we have here, namely instead of full amalgamation it was assumed that only unions of limit models are amalgamation bases. A follow up is in Grossberg, VanDieren, and Villaveces \cite{gvv}.

\begin{defin}
\begin{enumerate}

	\item $K$ has unique limit models in $\lambda$ iff if $M, N_1, N_2 \in K_\lambda$ and $\alpha_1, \alpha_2 < \lambda^+$ so that $N_\ell$ is $(\lambda, \alpha_\ell)$-limit over $M$, then $N_1 \cong_M N_2$.

\end{enumerate}
\end{defin}

It is an easy exercise to show that (2) holds if $\cf \alpha_1 = \cf \alpha_2$.  While many of the above papers prove the uniqueness of limit models in different contexts, the most relevant for our context is the proof that is outlined in \cite[Section II.4]{shelahaecbook} and detailed in Boney \cite{extendingframes}.  There, Shelah's frames are used to create a matrix of models to show that limit models are isomorphic.  Inspecting the proof, the only property used that is not a part of an independence property is a stronger continuity restricted to universal chains.  This follows from universal local character.

\begin{fact}[\cite{extendingframes}.8.2] \label{efthm}
If $\dnf$ is an independence relation so $\kappa^*_\omega(\dnf) = \delta$, then any two limit models of length at least $\delta$ are isomorphic.  Thus, if $\kappa^*_\omega(\dnf) = \omega$, then $K$ has unique limit models.
\end{fact}

\begin{remark}
It is unclear if this theorem is indeed an improvement of \cite{vandierennomax}, \cite{nomaxerrata} and  \cite{gvv}.  In Theorem \ref{efthm} full amalgamation is used.  Notice that the proof of Theorem \ref{efthm} uses that $\dnf$ is symmetric.  The other approaches for uniqueness of limit models don't use symmetry.  Additionally, the necessary existence and extension properties for nonsplitting can be proved from more natural hypotheses than we have here.

A proof for the original uniqueness statement from \cite{shvi635}  is still not known.
\end{remark}

\begin{cor} \label{ulmcor}
Suppose there is some $\kappa > LS(K)$ so that
\begin{enumerate}
	\item $K$ is fully $< \kappa$-tame and -type short;
	\item $K$ doesn't have the weak $\kappa$-order property; 
	\item $\dnf$ satisfies $(E)$; and
	\item it is categorical in some $\lambda > \kappa$
\end{enumerate}
Then $K$ has a unique limit model in each size in $[\kappa, \lambda)$.  Moreover, if $\lambda$ is a successor, then $K$ has unique limit models in each size above $\kappa$.
\end{cor}

{\bf Proof:} The first part follows from Theorems \ref{fulltheorem}, \ref{order2heir} and \ref{ulcthm} and Fact \ref{efthm}.  The moreover follows from the categoricity transfer of \cite[Theorem 5.2]{tamenessthree}. \hfill \dag

Note that the uniqueness of limit models as stated does not follow trivially from categoricity because it requires that the isomorphism fixes the base.

%%%%%%%%%%%%%%%%%%%%%%%%%%%%%%%%%%%%%%%%%%%%%%%
\section{The U-Rank} \label{uranksection}
%%%%%%%%%%%%%%%%%%%%%%%%%%%%%%%%%%%%%%%%%%%%%%%

Here we develop a $U$-rank for our forking and show that, under suitable conditions, it behaves as desired.  The $U$-rank was first introduced by Lascar \cite{lascar1975} for first order theories and first applied to AECs by \cite{sh394}.  They have also been studied by Hyttinen, Kesala, and Lessman in various non-elementary contexts; see \cite{lessman2000}, \cite{lessman2003}, \cite{hyttinenlessman2002}, and \cite{hyttinenkesala2006}.

For this section, we work with an abstract  independence relation $\adnf$ in the sense of Definition \ref{indepreldef}.  As a corollary, this means that, if the hypotheses of Theorem \ref{mainthm} hold, then the following holds for the $U$-rank defined in terms of $\kappa$-coheir.

\begin{hypothesis}
Suppose that $\adnf$ is an independence relation.  Recall that Hypothesis \ref{hyp} is also in effect.
\end{hypothesis}

Because we are working with $\adnf$ instead of $\dnf$, we will say that $gtp(A/N)$ does not $*$-fork over $M$ to mean $A \adnf_M N$.

\begin{defin}
We define $U$ with domain a type and range an ordinal or $\infty$ by, for any $p \in gS(M)$
\begin{enumerate}

	\item $U(p) \geq 0$;
	
	\item $U(p) \geq \alpha$ limit iff $U(p) \geq \beta$ for all $\beta < \alpha$;
	
	\item $U(p) \geq \beta + 1$ iff there is $M' \succ M$ with $\|M'\| = \|M\|$ and $p' \in gS(M')$ such that $p'$ is a $*$-forking extension of $p$ and $U(p') \geq \beta$;
	
	\item $U(p) = \alpha$ iff $U(p) \geq \alpha$ and $\neg( U(p) \geq \alpha+1)$; and
	
	\item $U(p) = \infty$ iff $U(p) \geq \alpha$ for every $\alpha$.

\end{enumerate}
\end{defin}

First we prove a few standard rank properties.  The first several results are true without the clause about the sizes of the model, but this is necessary later when we give a condition for the finiteness of the rank for Lemma \ref{ordinalbound}.

\begin{lemma}[Monotonicity] \label{rankmonotonicity}
If $M \prec N$, $p \in gS(M)$, $q \in gS(N)$, and $p \leq q$, then $U(q) \leq U(p)$.
\end{lemma}

{\bf Proof:}  We prove by induction on $\alpha$ that $p \leq q$ implies that $U(q) \geq \alpha$ implies $U(p) \geq \alpha$.  For limit $\alpha$, this is clear, so assume $\alpha = \beta + 1$ and $U(q) \geq \beta + 1$.  Then there is a $N' \succ N$ and $q^+ \in gS(N')$ that is a $*$-forking extension of $q$ and $U(q^+) \geq \beta$.  By $(M)$, it is also a $*$-forking extension of $p$.  Then $U(p) \geq \alpha$ as desired. \hfill \dag\\

\begin{lemma}[Invariance]\label{rankinvariance}
If $f \in Aut{ }\sea$ and $p \in gS(M)$, then $U(p) = U(f(p))$.
\end{lemma}

\begin{prop}[Ultrametric]
The $U$ rank satisfies the ultrametric property; that is, if we have $M \prec N_i$, $p \in gS(M)$ and distinct $\seq{q_i \in gS(N_i) \mid i < \alpha}$ are such that $a \models p$ iff there is an $i_0 < \alpha$ such that $a \models q_{i_0}$, then we have $U(p) = \max_{i < \alpha} U(q_i)$.
\end{prop}

Note that, as always, we assume $\alpha$ is well below the size of the monster model.

{\bf Proof:} We know that $p \leq q_i$ for all $i < \alpha$, so, by Lemma \ref{rankmonotonicity}, we have $\max_{i < \alpha} U(q_i) \leq U(p)$.  Since we have a monster model, we can find some $N^* \in K$ that contains all $N_i$.  By $(E)$, we can find some $p^+ \in gS(N^*)$ such that $p^+$ is a non-forking extension of $p$.  Now, let $a \models p^+$.  Since $p \leq p^+$, $a \models p$.  Since $p(\sea) = \cup_{i < \alpha} q_i(\sea)$, there is some $i_0 < \alpha$ such that $a \models q_{i_0}$.  But then $a \dnf_M N^*$ implies $a \dnf_M N_{i_0}$ by $(M)$, so $gtp(a/N_{i_0}) = q_{i_0}$ does not fork over $M$.  Then, by Theorem \ref{nf=U-rank} (proved independently), we have
$$U(p) = U(q_{i_0}) = \max_{i < \alpha} U(q_i)$$
\hfill \dag\\

We show that same rank extensions correspond exactly to non-forking  when the $U$-rank is ordinal valued.  One direction is clear from the definition.  For the other, we generalize first order proofs to the AEC context; our proof follows the one for \cite[Proposition 5.13]{pillay}.  First, we prove the following lemma.

\begin{lemma}\label{lemma7.6}
Let $N_0 \prec N_1 \prec \bar N_1$, $N_0 \prec \bar N_0 \prec \bar N_1$, and $N_0 \prec N_2$ be models with some $c \in \bar N_0$.  If
$$N_1 \adnf_{N_0} \bar N_0 \te{ and } N_2 \adnf_{\bar N_0} \bar N_1$$
then there is some $N_3$ extending $N_1$ and $N_2$ such that
$$c \adnf_{N_2} N_3$$
\end{lemma}

{\bf Proof:} We can use $(S)$ twice on $N_2 \adnf_{\bar N_0} \bar N_1$ to find $\bar N_2$ extending $N_2$ and $\bar N_0$ such that $\bar N_2 \adnf_{\bar N_0} \bar N_1$.  This contains $c$, so $(M)$ implies that $N_2 c \adnf_{\bar N_0}  N_1$.  By applying $(S)$ to the other non-$*$-forking from our hypothesis, we know $\bar N_0 \adnf_{N_0} N_1$.  By $(T)$, this means that $N_2 c \adnf_{N_0} N_1$.

Applying $(S)$ to this, there is some $N_3'$ extending $N_2$ and containing $c$ such that $N_1 \adnf_{N_0} N_3'$.  By $(M)$, we have that $N_1 \adnf_{N_2} N_3'$.  Applying $(S)$ one final time, we can find an $N_3$ extending $N_1$ and $N_2$ such that $c \adnf_{N_2} N_3$. \hfill \dag\\

\begin{theorem}\label{nf=U-rank}
Let $p \in gS(M_0)$ and $q \in gS(M_1)$ such that $p \leq q$ and $U(p), U(q) < \infty$.  Then
$$U(p) = U(q) \iff \te{$q$ is a non-$*$-forking extension of $p$}$$
\end{theorem}

{\bf Proof:} By definition, $U(p) = U(q)$ implies $q$ does not $*$-fork over $M_0$.  For the other direction, we show by induction on $\alpha$ that, for any $q$ that is a non-$*$-forking extension of some $p$, $U(p) \geq \alpha$ implies $U(q) \geq \alpha$. 

If $\alpha$ is $0$ or limit, this is straight from the definition.

Suppose that $U(p) \geq \alpha + 1$.  Then, there are $M_2 \succ M_0$ and $p_1 \in gS(M_2)$ such that $p_1$ is a $*$-forking extension of $p$ and $U(p_1) \geq \alpha$.

{\bf Claim:} We may pick $M_2$ and $p_1$ such that there is a $M_3$ extending $M_1$ and $M_2$ and $q_1 \in gS(M_3)$ so 
\begin{itemize}
	\item $q_1 \geq q, p_1$; and
	\item $q_1$ does not $*$-fork over $M_2$.
\end{itemize}

{\bf This Claim is enough:} Assume for contradiction that $q_1$ does not $*$-fork over $M_1$.  By \cite[Lemma 5.9]{bgkv}, a right version of transitivity also holds for $\adnf$:
\begin{center}
if $A \adnf_{M_0} M_1$ and $A \adnf_{M_1} M_2$ with $M_0 \prec M_1 \prec M_2$, then $A \adnf_{M_0} M_2$
\end{center}
Thus, $q_1$ would also not $*$-fork over $M_0$.  By $(M)$, this would imply that $p_1$ does not $*$-fork over $M_0$, a contradiction.  

Thus, $q_1$ is a $*$-forking extension of $q$.  Since $p_1 \leq q_1$ and $U(p_1) \geq \alpha$, Lemma \ref{rankmonotonicity} implies that $U(q_1) \geq \alpha$.  Thus, $U(q) \geq \alpha + 1$.

{\bf Proof of claim:} Let $d$ realize $q$ and $d'$ realize $p_1$.  Since both of these types extend $p$, there is some $f \in Aut_{M_0} \sea$ such that $f(d') = d$.  Set $M_2' = f(M_2)$.  We know that $d \adnf_{M_0} M_1$, so by $(S)$, there is some $\bar M_0 \succ M_0$ that contains $d$ so $M_1 \adnf_{M_0} \bar M_0$.  Pick $\bar M_1 \prec \sea$ that contains $\bar M_0$ and $M_1$.  By $(E)$, there is some $M_2''$ so that $gtp(M_2'/\bar M_0) = gtp(M_2''/\bar M_0)$ and $M_2'' \adnf_{\bar M_0} \bar M_1$.  Let $g \in Aut_{\bar M_0} \sea$ such that $g(M_2') = M_2''$; note that this fixes $d$.

We may now apply Lemma \ref{lemma7.6}.  This means there is some $M_3$ that extends $M_2''$ and $M_1$ such that $d \adnf_{M_2''} M_3$.  Now this proves our claim with $M_2''$ and $gtp(d/M_2'') = g(f(p_1))$ and witnesses $M_3$ and $q_1 = gtp(d/M_3)$. \hfill \dag\\

We now give a condition for the $U$ rank to be ordinal valued, as in \cite[Section 5]{sh394}.  First, note that clause about the model sizes in the definition of $U$ gives a bound for the rank.

\begin{lemma}[Ordinal Bound] \label{ordinalbound}
For each $\mu \geq \kappa$, there is some $\alpha_{K, \mu} < (2^\mu)^+$ such that for any $M \in K_\mu$ and $p \in gS(M)$, if $U(p) \geq \alpha_{K,\mu}$, then $U(p) = \infty$.
\end{lemma}

{\bf Proof:}   First, by Lemma \ref{rankinvariance}, the restriction to models of the same size means that there are at most $2^\mu$-many ordinals that are $U$-ranks of types over models of size $\mu$.  Second, there can be gaps below an ordinal valued rank; that is, given $N \in K_\mu$ and $q \in gS(N)$, if $\alpha < U(q) < \infty$, then there is a $*$-forking extension $q'$ of $q$ such that $U(q') = \alpha$.  If this were not the case, then a construction like the one in the next theorem would imply that $q$ has rank $\infty$.

Putting these together, we have that $\{U(p) : p \in gS(M), M\in K_\mu\} - \{\infty\}$ is an interval of size at most $2^\mu$.  Setting $\alpha_{K, \mu}$ to be the first ordinal not in this interval finishes the proof.\hfill \dag\\

This bound allows us to give a characterization of superstability in terms of an ordinal bound on the $U$-rank.

\begin{theorem}[Superstability]\label{U rank bounded iff nf wf}
Let $M \in K_\mu$ and $p \in gS(M)$.  Then the following are equivalent:
\begin{enumerate}

	\item $U(p) = \infty$.
	
	\item There is an increasing sequence of types $\seq{p_n : n < \omega}$ such that $p_0 = p$ and $p_{n+1}$ is a $*$-forking extension of $p_n$ for all $n < \omega$.

\end{enumerate}
\end{theorem}

{\bf Proof:} First, suppose $U(p) = \infty$ and set $p_0 = p$.  We will construct our sequence by induction such that $U(p_n) = \infty$.  Then $U(p_n) \geq \alpha_{K, \mu} + 1$, so there is a $*$-forking extension $p_{n+1}$ with the same sized domain and $U(p_{n+1}) \geq \alpha_{K,\mu}$.  But then $U(p_{n+1}) = \infty$ by Lemma \ref{ordinalbound} and our induction can continue.

Second, suppose we have such a sequence $\seq{p_n : n < \omega}$ and we will show, by induction on $\alpha$, that $U(p_n) \geq \alpha$ for all $n < \omega$.  The $0$ and limit stages are clear.  At stage $\alpha + 1$, $p_{n+1}$ is a $*$-forking extension of $p_n$ with rank at least $\alpha$.  Thus, $U(p_n) \geq \alpha + 1$. \hfill \dag\\

Ranks in a tame AEC have also been explored by Lieberman \cite{liebermanrank}.  Under a tameness assumption, he introduces a series of ranks that emulate Morley's Rank.

\begin{defin}[\cite{liebermanrank}.3.1]
Let $\lambda \geq \kappa$, where $K$ is $\kappa$-tame.  For $M \in K_\lambda$ and $p \in gS(M)$, we define $R^\lambda(p)$ inductively by
\begin{itemize}

	\item $R^\lambda[p] \geq 0$;
	
	\item $R^\lambda[p] \geq \alpha$ for limit $\alpha$ iff $R^\lambda[p] \geq \beta$ for all $\beta < \alpha$; and
	
	\item $R^\lambda[p] \geq \beta + 1$ iff there is $M' \succ M$ and distinct $\seq{p_i \in gS(M'): i < \lambda^+}$ such that $p \leq p_i$ and $R^\lambda[p_i] \geq \beta$ for all $i < \lambda^+$.

\end{itemize}

If $\|M\| > \lambda$ and $p \in gS(M)$, then
$$R^\lambda[p] = \min \{ R^\lambda[p \rest N] : N \prec M, \|N\| = \lambda \}$$
\end{defin}

The $U$-rank for any independence relation dominates these Morley Ranks at least for domains of size $\lambda$.  Thus, the finiteness of the $U$-rank, which follows from local character, implies that an AEC is totally transcendental and that the stability transfer results of \cite[Section 5]{liebermanrank} apply.

\begin{theorem}
Suppose $K$ is $<\kappa$-tame.  Let $M \in K_\lambda$, $p \in gS(M)$, and $\lambda \geq \kappa$.  Then $U(p)\geq R^\lambda(p)$.
\end{theorem}

{\bf Proof:} We prove, simultaneously for all types, that $R^\lambda(p) \geq \alpha$ implies $U(p) \geq \alpha$ for all $\alpha$ by induction.  For $\alpha = 0$ or limit, this is easy.

Suppose $R^\lambda(p) \geq \alpha+1$.  Let $M'$ and $\seq{p_i \in gS(M') : i < \lambda^+}$ witness this.  $p$ has a unique $*$-nonforking extension to $M'$, call it $p^*$.  Thus, almost all of the $p_i$ $*$-fork over $M$; let $p_{i_0}$ be one of them.  Then, $p_{i_0} \neq p^*$, so by $<\kappa$-tameness, there is some $M_0 \prec M'$ of size $< \kappa$ such that $p_{i_0} \rest M_0 \neq p^* \rest M_0$.  Let $M'' \prec M'$ contain $M$ and $M_0$ such that $\|M\| = \|M''\|$ and $p' = p_{i_0} \rest M''$.  Then
\begin{itemize}

	\item $p'$ extends $p$;
	
	\item $p'$ is a $*$-forking extension of $p$ because it differs from the non-$*$-forking extension, $p^* \rest M_{i_0}$; and
	
	\item $R^\lambda(p') \geq R^\lambda(p_{i_0})$ by \cite[Proposition 3.3]{liebermanrank}.  So $R^\lambda(p') \geq \alpha$.  By induction, this means $U(p') \geq \alpha$.

\end{itemize}
So $U(p) \geq \alpha+1$, as desired.  \hfill \dag\\

%%%%%%%%%%%%%%%%%%%%%%%%%%%%%%%%%%%%%%%%%%%%%%%
\section{Large cardinals revisited} \label{lcrevisitedsection}
%%%%%%%%%%%%%%%%%%%%%%%%%%%%%%%%%%%%%%%%%%%%%%%

In this section, we discuss the behavior of non-forking in the presence of large cardinals.  We return to just assuming Hypothesis \ref{hyp}, that $K$ satisfies amalgamation, joint embedding, and no maximal models.

Recall that $\kappa$ is strongly compact iff $\kappa$ is regular and every $\kappa$ complete filter can be extended to a $\kappa$ complete ultrafilter; see \cite[Chapter 20]{jech} for a reference.  Boney \cite[Theorem 4.5]{tamelc} proved that the tameness and type shortness hypotheses of Theorem \ref{mainthm} follow from $\kappa$ being a strongly compact cardinal.

\begin{fact}[\cite{tamelc}] \label{tamelcsc} \label{strongcompactnesstheorem}
If $\kappa$ is strongly compact and $K$ is an AEC with $LS(K) < \kappa$, then $K$ is fully $<\kappa$ tame and fully $<\kappa$ type short.
\end{fact}

A similar result holds for AECs axiomatized in $L_{\kappa, \omega}$.  However, \cite{makkaishelah} deals with this case more fully, so we focus on $LS(K) < \kappa$.  In fact, extending the results of \cite{makkaishelah} to general AECs via the methods of \cite{tamelc} was the motivation for this paper.  The key tool of \cite{tamelc} is a \L o\'{s}' Theorem for AECs (see \cite[Theorems 4.3 and 4.7]{tamelc}) that says that such AECs are closed under sufficiently complete ultraproducts, that ultraproducts of embeddings is an embedding of the ultraproducts, and more.

We now detail a construction that will be used often in the following proof.  This construction and the proof of the following theorem draw inspiration from \cite[Proposition 4.9]{makkaishelah}.\\

Suppose that $M \prec N$ and $U$ is a $\kappa$ complete ultrafilter over $I$.  Then \L o\'{s}' Theorem for AECs states that the canonical ultrapower embedding $h: N \to \Pi N/U$ that takes $n$ to the constant function $[i \mapsto n]_U$ is a $K$-embedding.  We can expand $h$ to some $h^+$ that is an $L(K)$-isomorphism with range $\Pi N/U$ and set $N^U := (h^+)^{-1}[\Pi N/U]$.  This is a copy of the ultraproduct that actually contains $N$.  Similarly, we can set $M^U := (h^+)^{-1}[\Pi M/U]$.  The following claim is key.

{\bf Claim:} $M^U \dnf_M N$.

{\bf Proof:} Let $N^- \prec N$ be small and $a \in M^U$.  Then $h^+(a) = [f]_U$ for some $[f]_U \in \Pi M/U$.  Denote $gtp(a/N^-)$ by $p$.  Then, by \cite[Theorem 4.7]{tamelc}, we have
\begin{eqnarray*}
a &\vDash& p\\
h^+(a) = [f]_U &\vDash& h^+(p) = h(p)\\
X&:=& \{i \in I : f(i) \vDash p \} \in U
\end{eqnarray*}
Since $[f]_U \in \Pi M/U$, there is some $i_0 \in X$ such that $f(i_0) \in M$.  Then $f(i_0) \vDash p$ as desired. \hfill \dag\\

We now show that non-forking is very well behaved in the presence of a strongly compact cardinal.  Note that the second part says that the local character property follows from Existence and the third part improves on Theorem \ref{catexist} by showing that categoricity implies an analogue of superstability instead of just an analogue of simplicity.

\begin{theorem}\label{sceasy}
Suppose $\kappa$ is strongly compact and $K$ is an AEC such that $LS(K) < \kappa$.  Then
\begin{enumerate}

	\item $\dnf$ satisfies Extension.
	
	\item Suppose $\dnf$ satisfies Existence.  Let $M = \cup_{i < \alpha} M_i$ such that $\alpha = \cf \alpha$.  If
	\begin{itemize}
		\item $|A| < \alpha$ when $|A| < \kappa$; or 
		\item $|A|^{<\kappa} < \alpha$ when $|A| \geq \kappa$,
	\end{itemize}	
	then there is $i_* < \alpha$ such that $A \dnf_{M_{i_*}} M$.
	
	\item If $K$ is categorical in some $\lambda = \lambda^{<\kappa}$, then $\dnf$ is an independence relation with $$\kappa_\alpha(\dnf) \leq \begin{cases} \omega+|\alpha| & \alpha < \kappa\\ |\alpha|^{<\kappa} & \alpha \geq \kappa \end{cases}$$
	Moreover, Hypothesis \ref{hyp} can be weakened to just referring to $K_{<\kappa}$.

\end{enumerate}
\end{theorem}

{\bf Proof:}  
\begin{enumerate}
	\item Suppose that $A \dnf_{M_0} N$ and let $N^+ \succ N$.  This means that every $<\kappa$-approximation to $gtp(A/N)$ is realized in $M_0$.  We use the same argument as \cite[Theorem 4.10]{tamelc}, even though this is not a type over $M_0$.  Let $U$ be a $\kappa$-complete, fine ultrafilter over $P_\kappa A \times P^*_\kappa N$ and let $\seq{b_a^{A_0, N^-} \in M_0 : a \in A_0}$ realize $gtp(A_0/N^-)$.  Then $A' = \seq{[(A_0, N^-) \mapsto b_a^{A_0, N^-}]_U : a \in A}$ is a realization of $h(gtp(A/N))$ by \cite[Theorem 4.7]{tamelc}.  By the above claim, $M_0^U \dnf_{M_0} N^+$.  By $(M)$, this implies $h^{-1}(A') \dnf_{M_0} N^+$.  Since $gtp(A/N) = gtp(h^{-1}(A')/N)$, this finishes the proof.
	
	\item Note that $A \dnf_M M$ by Existence.  First suppose that $|A| < \kappa$.  We break into cases based on the size of $\alpha$.
	
	If $\alpha < \kappa$, then, as before, we can use the fact that $A \dnf_M M$ to find a $\kappa$-complete ultrafilter $U$ on $I$ such that $p$ is realized in $M^U$.  Since $\alpha < \kappa$ and $U$ is $\kappa$-complete, we have that $M^U = \cup_{i < \alpha} M_i^U$.  Let $A' \in M^U$ realize $p$.  Since $|A'| <  \alpha$ and $\alpha$ is regular, there is some $i_* < \alpha$ such that $A' \in M_{i_*}^U$.  Thus, by the claim,
	$$M^U_{i_*} \dnf_{M_{i_*}} M$$
	$$A' \dnf_{M_{i_*}} M$$
	By $(I)$, $A \dnf_{M_{i_*}} M$.

	Now suppose that $\alpha \geq \kappa$.  For contradiction, suppose that $ A \not \dnf_{M_i} M$ for all $i < \alpha$.  We now build an increasing and continuous sequence of ordinals $\seq{i_j : j < |A|^+ }$ by induction.  Let $i_0 = 0$.  Given $i_j$, we know that $p$ forks over $M_{i_j}$.  By the definition, there is a small $M^- \prec M$ such that $gtp(A/M^-)$ is not realized in $M_{i_j}$.  Since $ \alpha \geq \kappa$ is regular, there is some $i_{j+1} > i_j$ such that $M^- \prec M_{i_{j+1}}$.  Then $A \not\dnf_{M_{i_j}} M_{i_{j+1}}$.  Set $M^* = \cup_{j < |A|^+} M_{i_j}$.  Then, by Monotonicity, $p \rest M^*$ forks over $M_{i_j}$ for all $j < |A|^+$.  Since $|A|^+ < \kappa$, this contradicts the previous paragraph.
	
	 Second suppose that $|A| \geq \kappa$.  For each $a \in {}^{<\kappa} A$, there is some $i_a < \alpha$ such that $a \dnf_{M_{i_a}} M$.  Set $i_* = \sup i_a$; by assumption $i_* < \alpha$.  Then $A \dnf_{M_{i_*}} M$. 
	
	\item From inaccessibility, we know that $\sup_{\mu<\kappa}(\beth_{(2^{\mu})^+}) = \kappa$, so Existence holds by Theorem \ref{catexist}.  Then Extension holds by the first part, so $(E)$ holds.  Theorem \ref{strongcompactnesstheorem}  tells us that $K$ is $< \kappa$-tame and -type short.  Finally, as outlined in the discussion after Theorem 5.1, the weak $\kappa$-order property with $\kappa$ inaccessible implies many models in all cardinals above $\kappa$, which is contradicted by categoricity in $\lambda$.\\
	
	For the moreover, Baldwin and Boney \cite[Corollary 3.16]{hanfap} show that amalgamation and joint embedding on $K_{<\kappa}$ imply those properties for all of $K$.  No maximal models above $\kappa$ follows from strong compactness by taking nonprincipal ultrapowers.\hfill \dag\\
\end{enumerate}

Additionally, with the full strength of a strongly compact cardinal, we can reprove much or all of \cite[Section 4]{makkaishelah} in an AEC context.  One complication is that \cite[Definition 4.23]{makkaishelah} defines weakly orthogonal types by having an element in the nonforking relation where we require a model.  However, this definition has already been generalized at \cite[Section III.6]{shelahaecbook}.

\cite{tamelc} also proves weaker versions of Theorem \ref{strongcompactnesstheorem} from assumptions of measurable or weakly compact cardinals.  These in turn could be used to produce weaker versions of Theorem \ref{sceasy}.  However, \cite{makkaishelah} is not the only time independence relations have been studied in infinitary contexts with large cardinals.  Kolman and Shelah \cite{kosh} and Shelah \cite{measure2} investigate the consequences of categoricity in $L_{\kappa, \omega}$ when $\kappa$ is measurable.  In \cite{kosh}, they use $\kappa$-complete ultralimits.  They denote such an ultralimit of $M$ by $Op(M)$ and the canonical embedding by $f_{Op}:M \to Op(M)$.  In \cite{measure2}, Shelah introduces the following independence relation. 

\begin{defin}[\cite{measure2}.1.5]
Let $K$ be essentially below $\kappa$ measurable.  Define a 4-place relation $\dnks$ by $M_1 \dnks_{M_0}^{M_3} M_2$ iff there is an ultralimit operation $Op$ with embedding $f_{Op}$ and $h:M_3 \to Op(M_1)$ such that the following commutes

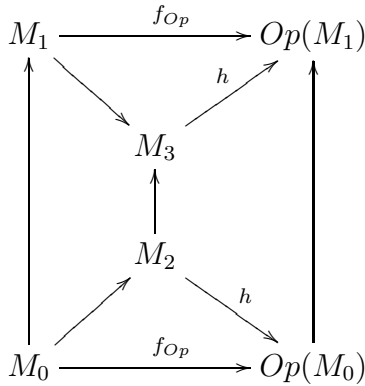
\[
\xymatrix{
M_1 \ar[rr]^{f_{Op}} \ar[dr] & & Op(M_1) \\
& M_3 \ar[ur]^h & \\
& M_2 \ar[u] \ar[dr]^h & \\
M_0 \ar[uuu] \ar[ur] \ar[rr]^{f_{Op}} & & Op(M_0) \ar[uuu]
}
\]

\end{defin}
We would like to show that $\kappa$-coheir and Shelah's nonforking are dual\footnote{By Proposition \ref{heircoduality}, this would mean it is equivalent to heir.} in all situations.  However, this does not seem to be the case as the proofs depend on the existence of regular ultrafilters or structure results from \cite{measure2}.  One direction always holds.

For the following results, we assume that $K$ is an AEC essentially below\footnote{This means $LS(K) < \kappa$ or $K$ is axiomatized by a $L_{\kappa, \omega}$-theory and ordered by elementary submodel according a suitable fragment; see \cite[Definition 2.10]{tamelc}.} $\kappa$.

\begin{prop}
If $\kappa$ is measurable, then, for all $M_0, M_1, M_2 \in K$, 
$$ \exists M_3\te{ so } M_2 \dnks_{M_0}^{M_3} M_1  \implies M_1 \dnf_{M_0} M_2 $$
\end{prop}

{\bf Proof:} Suppose that there is an $M_3$ such that $M_2 \dnks_{M_0}^{M_3} M_1$.  The claim above generalized to ultralimits implies $f_{Op}^{-1}(Op(M_0)) \dnf_{M_0} M_2$ for the ultralimit witnessing $M_2$.  We have that $h:M_1 \to Op(M_0)$, so, by Monotonicity, we have that $f_{Op}^{-1}(h(M_1)) \dnf_{M_0} M_2$.  By the diagram, $f_{Op}^{-1} \circ h$ fixes $M_2$, we have that $gtp(f_{Op}^{-1}(h(M_1))/M_2) = gtp(M_1/M_2)$.  By Invariance, this means that $M_1 \dnf_{M_0} M_2$. \hfill \dag\\

For the other direction, if $\kappa$ is simply measurable, then the result only holds on $K_\kappa$.

\begin{prop}
If $\kappa$ is measurable, then, for all $M_0, M_1, M_2 \in K_\kappa$, 
$$M_1 \dnf_{M_0} M_2  \implies \exists M_3\te{ so } M_2 \dnks_{M_0}^{M_3} M_1$$
\end{prop}

{\bf Proof:} Suppose that $M_1 \dnf_{M_0} M_2$.  Adapting the proof that every countably incomplete measure is $\omega$-regular (see \cite[Proposition 4.3.4]{changkeisler}), every maximally $\kappa$-complete ultrafilter is $\kappa$-regular.  Let $U$ be a maximally $\kappa$-complete ultrafilter (such as the one arising from a measurable embedding)% and $\{X_\alpha \in U : \alpha < \kappa\}$ be a regularizing family
.  Adapting first order arguments (see \cite[Theorem VI.1.4]{shelahfobook}), we have the following.

{\bf Claim:} For $M \in K_\kappa$, if $N \in K_\kappa$ and $p \in gS^\kappa(N)$ is a Galois type such that every small subtype is realized in $M$, then $p$ is realized in $M^U$.

%To see the claim, find a resolution $\seq{N_i \in K_{<\kappa} : i < \kappa}$ of $N$.  For each $\alpha < \kappa$, there is $\seq{a^\alpha_\beta \in M : \beta < \alpha}$ that realizes $p_\alpha := p^\alpha \rest N_\alpha$ by assumption.  For every $\alpha < \kappa$, set $Y_\alpha = \{ i < \kappa : \alpha \in X_i\}$; by regularity, this has size $< \kappa$.  For $i < \kappa$, define $f_i \in \Pi M$ by XXXX

By the claim at the start of the section, we have that $M^U_0 \dnf_{M_0} M_2$.  Since $M_1 \dnf_{M_0} M_2$, the hypothesis of the above claim is satisfied with $M_0$ and $gtp(M_1/M_2)$.  Thus, it is realized in $M_0^U$ and there is $f \in Aut_{M_2} \sea$ such that $f(M_1) \prec M_0^U$.  Set $M_3 = f^{-1}[M_2^U]$.  Then we have the following commuting diagram:

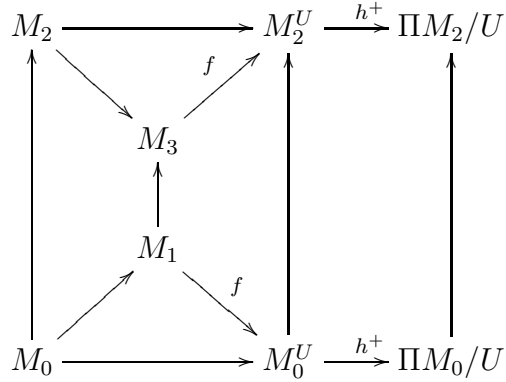
\[
\xymatrix{
M_2 \ar[rr] \ar[dr] & & M_2^U \ar[r]^{h^+} & \Pi M_2/U \\
& M_3 \ar[ur]^f & &\\
& M_1 \ar[u] \ar[dr]^f & &\\
M_0 \ar[uuu] \ar[ur] \ar[rr] & & M_0^U \ar[uuu] \ar[r]^{h^+}& \Pi M_0/U \ar[uuu]
}
\]

Collapsing this diagram gives

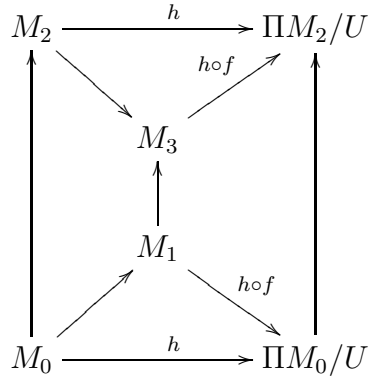
\[
\xymatrix{
M_2 \ar[rr]^h \ar[dr] &  & \Pi M_2/U \\
& M_3 \ar[ur]^{h \circ f} & \\
& M_1 \ar[u] \ar[dr]^{h \circ f} & \\
M_0 \ar[uuu] \ar[ur] \ar[rr]^h & & \Pi M_0/U \ar[uuu]
}
\]
Note that an ultraproduct is a suitable ultralimit operation and the ultrapower embedding is its corresponding embedding.  Thus $M_2 \dnks_{M_0}^{M_3} M_1$. \hfill\dag\\

As can be seen from the proof, moving beyond $K_\kappa$ would require $\kappa$-complete ultrafilters of higher regularity.  This seems to require more than just a measurable cardinal and, in particular, strongly compact is enough.

\begin{prop}
If $\kappa$ is $\lambda$-strongly compact, then, for all $M_0, M_1, M_2 \in K_{[\kappa, \lambda]}$, 
$$M_1 \dnf_{M_0} M_2  \implies \exists M_3\te{ so } M_2 \dnks_{M_0}^{M_3} M_1$$
\end{prop}

The proof of this proposition is nearly identical to the one above; the extra step is that the $\kappa$ complete, fine ultrafilter on $P_\kappa \lambda$ gives a $\lambda$ regular ultrafilter.  Another option would be to take advantage of continuity properties of $\dnks$ and use direct limits of $\kappa$ sized models; however, \cite{measure2} indicates that these continuity properties require additional hypotheses.

This gives the following corollary.

\begin{cor}
Let $K$ be an AEC essentially below $\kappa$ strongly compact and let $M_0 \prec M_1, M_2 \in K$.  Then
$$M_1 \dnf_{M_0} M_2 \iff \exists M_3\te{ so } M_2 \dnks_{M_0}^{M_3} M_1$$
\end{cor}

%In Section \ref{lcsection}, we saw that knowledge of the relationship between two possible independence relations (nonforking and heiring) enabled us to find information about the local character of nonforking.  This raises the possibility that studying what implications hold between the various candidates--nonforking, nonsplitting, etc--will lead to new information about these relations.  This will be done in \cite{relationrelation}.

%\bibliography{}{}

\begin{thebibliography}{She01b}


\bibitem[Ad09]{Adler}  Hans Adler. A geometric introduction to forking and thorn-forking, Journal of Mathematical Logic {\bf 9} (2009) (1), 1--20.


\bibitem[Bal09]{baldwinbook}
John Baldwin, \emph{{\bf Categoricity}}, University Lecture Series, American
  Mathematical Society, 2009.

\bibitem[BaBo]{hanfap}
John Baldwin and Will Boney, \emph{Hanf numbers and presentation theorems in AECs}, In preparation.

\bibitem[BET07]{nperp}
John Baldwin, Paul Eklof, and Jan Trlifaj, \emph{${}^{\perp}N$ as an AEC}. Annals of Pure and Applied Logic {\bf 149} (2007)

\bibitem[BK09]{untame}
John Baldwin and Alexei Kolesnikov, \emph{Categoricity, amalgamation, and
  tameness}, Israel Journal of Mathematics \textbf{170} (2009), no.~1,
  411--443.
  
%  \bibitem[BLS1003]{bls1003}
%John Baldwin, Paul Larson, and Saharon Shelah, \emph{Almost Galois $\omega$-stable classes}, To appear, Journal of Symbolic Logic, \verb+http://www.users.muohio.edu/larsonpb/BlLrSh1003sep28.pdf+

\bibitem[BlSh862]{nonlocality}
John Baldwin and Saharon Shelah, \emph{Examples of non-locality}, Journal of
  Symbolic Logic \textbf{73} (2008), 765--782.
  
\bibitem[BYCh14]{Ch1} Itai Ben-Yaacov and Artem Chernikov, \emph{An independence theorem for NTP2 theories}, Journal of Symbolic Logic {\bf 79} (2014), no. 1, 135-153.

\bibitem[Bona]{longtypes}
Will Boney, \emph{Computing the number of types of infinite length}, Accepted, Notre Dame Journal of Formal Logic, \verb+http://arxiv.org/abs/1309.4485+

\bibitem[Bonb]{gamult}
Will Boney, \emph{The $\Gamma$-ultraproduct and averageable classes}, Preprint, \verb+http://arxiv.org/abs/1511.00982+.

%\bibitem[Bonc]{tormod}
%Will Boney, \emph{Some model theory of torsion modules over PIDs}, In preparation.

\bibitem[Bon14a]{extendingframes}
Will Boney, \emph{Tameness and extending frames}, Journal of Mathematical Logic {\bf 14}, no. 2

\bibitem[Bon14b]{tamelc}
Will Boney \emph{Tameness from large cardinals axioms}, Journal of Symbolic Logic {\bf 79}, no. 4, 1092-1119.

\bibitem[BGKV]{bgkv}
Will Boney, Rami Grossberg, Alexei Kolesnikov, and Sebastien Vasey, \emph{Canonical Forking in AECs}, Annals of Pure and Applied Logic {\bf 167} (2016), no. 7, 590-613

\bibitem[BGVV]{svthm}
Will Boney, Rami Grossberg, Monica VanDieren, and Sebastien Vasey, \emph{Superstability from Categoricity in Abstract Elementary Classes}, In preparation.

\bibitem[BoUn]{bonung}
Will Boney and Spencer Unger, \emph{Large cardinals from tameness}, In preparation.

\bibitem[BoVa]{unionsat}
Will Boney and Sebastien Vasey, \emph{Chains of Saturated Models in AECs}, Submitted, \verb+http://arxiv.org/abs/1503.08781+.

  \bibitem[BuLe03]{BuLe}
Steve Buechler and Olivier Lessmann. Simple homogeneous models.
\emph{J. Amer. Math. Soc.}, {\bf 16} (2003), no. 1, 91--121 

\bibitem[ChKe90]{changkeisler} C. C. Chang and H. J. Keisler, \textbf{Model Theory}, Third Edition,
\newblock 208 pages, {North-Holland 1990}


\bibitem[CKS1007]{Ch2} Artem Chernikov, Itay Kaplan and Saharon Shelah. \emph{On non-forking spectra}, \verb+http://arxiv.org/pdf/1205.3101.pdf+, Submitted.

\bibitem[ChKa12]{ChKa} Artem Chernikov and Itay Kaplan.  \emph{Forking and dividing in NTP2 theories}, Journal of Symbolic Logic {\bf 77} (2012), 1-20.

\bibitem[Gro02]{grossberg2002}
Rami Grossberg, \emph{Classification theory for abstract elementary classes},
  Logic and Algebra (Yi~Zhang, ed.), vol. 302, American Mathematical Society,
  2002, pp.~165--204.

\bibitem[Gro1X]{ramibook}
\bysame, \emph{{\bf A Course in Model Theory}}, In Preparation, 201X.

\bibitem[GIL02]{primersimple}
Rami Grossberg, Jose Iovino, and Olivier Lessman, \emph{A primer of simple
  theories}, Archive of Mathematical Logic \textbf{41} (2002), 541--580.

\bibitem[GrLe02] {GrLe1} Rami Grossberg and Olivier Lessmann, \emph{
Shelah's Stability Spectrum and Homogeneity Spectrum},
 {  Archive for 
mathematical Logic}, {\bf 41}, (2002) 1, 1-31.

\bibitem[GrLe05]{GrLe05} Rami Grossberg and Olivier Lessmann,  \emph{Abstract decomposition theorem and applications}, Contemporary Mathematics {\bf 380} (2005), AMS, pp. 73--108.

\bibitem[GrSh]{GrSh}
Rami Grossberg and Saharon Shelah,
\newblock \emph{On universal locally finite groups}
\newblock { Israel J. of Math. }, {\textbf 44}, 1983, 289--302.

\bibitem[GV06a]{tamenessthree}
Rami Grossberg and Monica VanDieren, \emph{Categoricity from one successor
  cardinal in tame abstract elementary classes}, Journal of Mathematical Logic
  \textbf{6} (2006), no.~2, 181--201.

\bibitem[GV06b]{tamenessone}
Rami Grossberg and Monica VanDieren, \emph{Galois-stability for tame abstract elementary classes}, Journal
  of Mathematical Logic \textbf{6} (2006), no.~1, 25--49.

\bibitem[GV06c]{tamenesstwo}
Rami Grossberg and Monica VanDieren, \emph{Shelah's categoricity conjecture from a successor for tame
  abstract elementary classes}, Journal of Symbolic Logic \textbf{71} (2006),
  no.~2, 553--568.

\bibitem[GVV]{gvv}
Rami Grossberg, Monica VanDieren, and Andres Villaveces, \emph{Uniqueness of
  limit models in abstract elementary classes}, Submitted. \verb+http://www.math.cmu.edu/~rami/GVV_1-24b-15-no-vp.pdf+

\bibitem[GVas]{gvas}
Rami Grossberg and Sebastien Vasey.
\emph{Superstability in abstract elementary classes}, In preparation.


\bibitem[HaHa84]{hh84}
Victor Harnik and Leo Harrington, \emph{Fundamentals of Forcing}. Annals of Pure and Applied Logic {\bf 26} (1984), 245-286.

\bibitem[HaSh323]{hash323}
Bradd Hart and Saharon Shelah, \emph{Categoricity over {$P$} for first order
  {$T$} or categoricity for {$\phi\in{\rm L}_ {\omega_ 1\omega}$} can stop at 
  $\aleph_ k$ while holding for $\aleph_ 0,\cdots,\aleph_ {k-1}$}, Israel
  Journal of Mathematics \textbf{70} (1990), 219--235.

%\bibitem[HySh629]{strsplithom}
%Tapani Hyttinen and Saharon Shelah, \emph{Strong splitting in stable
%  homogeneous models}, Annals of Pure and Applied Logic \textbf{103} (2000),
%  201--228.
  
\bibitem[HyKe06]{hyttinenkesala2006}
Tapani Hyttinen and Meeri Kes{\''a}l{\''a}, \emph{Independence in finitary abstract elementary classes}, Annals of Pure and Applied Logic \textbf{143} (2006), 103-138
  
\bibitem[HyKe11]{hyttinenkesala11}
Tapani Hyttinen and Meeri Kes{\''a}l{\''a} (2011), \emph{Categoricity Transfer in Simple Finitary Abstract Elementary Classes}. Journal of Symbolic Logic {\bf 76} (2011), 759 - 806.

  \bibitem[HyLe02]{hyttinenlessman2002}
  Tapani Hyttinen and Olivier Lessman, \emph{A rank for the class of elementary submodels of a superstable homogeneous model}, The Journal of Symbolic Logic \textbf{67} (2002), 1469-1482

\bibitem[JrSh875]{jrsh875}
Adi Jarden and Saharon Shelah, \emph{Non-forking frames in abstract elementary classes}, Annals of Pure and Applied Logic \textbf{164} (2013), 135-191

\bibitem[Jec06]{jech}
Thomas Jech, \emph{{\bf Set Theory}}, 3rd ed., Springer, 2006.

\bibitem[Kan08]{kanamori}
Akihiro Kanamori, \emph{{\bf The Higher Infinite: Large Cardinals in Set Theory
  from Their Beginnings}}, 2nd ed., Springer, 2008.
  
  \bibitem[Kei71]{KeLomega} H. Jerome Keisler, \textbf{Model Theory for
Infinitary Logic}.
\newblock 208 pages, {North-Holland 1971}

  
  \bibitem[Kim98]{Kim}
Byunghan Kim.
\newblock \emph{Forking in simple unstable theories},
\newblock { J. of London Math. Society},
\newblock  (2) 57 (1998) 257--267.


\bibitem[KiPi97]{KP1}
Byunghan Kim and Anand Pillay,
\newblock \emph{Simple theories}, 
{ Annals of Pure and Applied Logic},
\newblock  88 (1997) 149-164.

\bibitem[KoSh362]{kosh}
Oren Kolman and Saharon Shelah, \emph{Categoricity of theories in {$L_{\kappa \omega}$}, when $\kappa$ is a measureable cardinal, part 1}, Fundamenta Mathematica \textbf{151} (1996), 209-240

\bibitem[Las75]{lascar1975}
Daniel Lascar, \emph{D\'{e}finissabilit\'{e} de types on th\'{e}orie des
  mod\`{e}les}, Ph.D. thesis, Universite Paris VIII, 1975.

\bibitem[Las76]{las76}
Daniel Lascar, \emph{Ranks and definability in superstable theories}, Israel Journal of Mathematics {\bf 23} (1976), no. 1, 53-87

\bibitem[Les00]{lessman2000}
Olivier Lessman, \emph{Ranks and pregeometries in finite diagrams}, Annals of Pure and Applied Logic \textbf{106} (2000), 49-83

%\bibitem[Les00b]{Le1} Olivier Lessmann. 
%Pregeometries in finite diagrams.
%\relax {\em Ann. Pure Appl. Logic}, {\bf 106},  no. 1-3, 49--83, 2000.

\bibitem[Les03]{lessman2003}
Olivier Lessman, \emph{Categoricity and U-rank in excellent classes}, The Journal of Symbolic Logic \textbf{68} (2003), 1317-1336

\bibitem[Lie13]{liebermanrank}
Michael Lieberman, \emph{Rank functions and partial stability spectra for tame
  abstract elementary classes}, Notre Dame Journal of Formal Logic \textbf{54}
  (2013), no.~2, 153--166.

\bibitem[MaSh285]{makkaishelah}
 Michael Makkai and Saharon Shelah, \emph{Categoricity of theories in
  {$L_{\kappa \omega}$}, with $\kappa$ a compact cardinal}, Annals of Pure and
  Applied Logic \textbf{47} (1990), 41--97.
  
\bibitem[Mar75]{marcus75}
Leo Marcus, \emph{A theory with only non-homogeneous models}, Algebra Universalis {\bf 5} (1975)

\bibitem[Pil83]{pillay}
Anand Pillay, \emph{{\bf An Introduction to Stability Theory}}, Dover
  Publications, Inc, 1983.

\bibitem[Sh:c]{shelahfobook}
Saharon Shelah, \emph{{\bf Classification theory and the number of nonisomorphic
  models}}, 2nd ed., vol.~92, North-Holland Publishing Co., Amsterdam,
  xxxiv+705 pp, 1990.


\bibitem[Sh:e]{shelahnonstructurebook}
Saharon Shelah, \emph{{\bf Non structure theory}}, In preparation, \verb+http://shelah.logic.at/nonstructure/+

\bibitem[Sh:g]{shg}
Saharon Shelah, \emph{{\bf Cardinal Arithmetic}}, Oxford Logic Guides no. 29, Oxford University Press, 1994.

\bibitem[Sh:h]{shelahaecbook}
Saharon Shelah, \emph{{\bf Classification Theory for Abstract Elementary Classes}},
  vol. 1 \& 2, Mathematical Logic and Foundations, no. 18 \& 20, College
  Publications, 2009.

\bibitem[Sh3]{Sh3}
Saharon Shelah.
\newblock Finite diagrams stable in power.
\newblock {\em Ann. Math. Logic}, {\bf 2}, 69--118, 1970/1971.

\bibitem[Sh31]{sh31}
Saharon Shelah, \emph{Categoricity of uncountable theories}, Proceedings of the Tarski
  Symposium, 1971, 1974, pp.~187--203.

\bibitem[Sh88]{sh88}
Saharon Shelah, \emph{Classification of nonelementary classes, {II}. {A}bstract
  elementary classes}, Classification theory (John Baldwin, ed.), 1987,
  pp.~419--497.

\bibitem[Sh300]{sh300}
Saharon Shelah, \emph{Universal classes}, Classification theory (John Baldwin, ed.),
  1987, pp.~264--418.

\bibitem[Sh394]{sh394}
Saharon Shelah, \emph{Categoricity for abstract classes with amalgamation}, Annals of
  Pure and Applied Logic \textbf{98} (1990), 261--294.



\bibitem[Sh472]{measure2}
Saharon Shelah, \emph{Categoricity of theories in {$L_{\kappa \omega}$}, when $\kappa$
  is a measureable cardinal, part 2}, Fundamenta Mathematica \textbf{170}
  (2001), 165--196.


\bibitem[Sh576]{sh576}
Saharon Shelah, \emph{Categoricity of an abstract elementary class in two successive
  cardinals}, Israel Journal of Mathematics (2001), no.~126, 29--128.

\bibitem[Sh600]{sh600}
Saharon Shelah, \emph{Categoricity in abstract elementary classes: going up
  inductively}, math.LO/0011215. See Chapter II of volume I [Sh:h].

  
  \bibitem[Sh705]{Sh705}
Saharon Shelah,
\newblock 
Toward classification theory of 
good $\lambda$ frames and abstract elementary classes.   See Chapter III of volume I [Sh:h].

  \bibitem[Sh734]{Sh734}
Saharon Shelah,
\newblock 
Categoricity and solvability of A.E.C., quite highly, 143 pages, \verb+http://arxiv.org/pdf/0808.3023.pdf+

\bibitem[Sh1019] {sh1019} Saharon Shelah, \emph{Model-theory for theta-complete ultrapowers}, \verb+http://arxiv.org/pdf/1303.5247v1.pdf+

\bibitem[ShVi635]{shvi635}
Saharon Shelah and Andres Villaveces, \emph{Toward categoricity for classes
  with no maximal models}, Annals of Pure and Applied Logic \textbf{97} (1999),
  1--25.


  
  \bibitem[Si]{Simon}
Pierre Simon. \emph{{\bf Lecture notes on NIP theories}} In preparation, \verb+http://www.normalesup.org/~simon/NIP_lecture_notes.pdf+


\bibitem[Vasa]{vaseytameframe}
Sebastien Vasey, \emph{Forking and superstability in tame AECs}, The Journal of Symbolic Logic {\bf 81} (2016), no. 1, 357?383

\bibitem[Vasb]{vaseybehemoth} 
Sebastien Vasey,
 \emph{Infinitary stability theory},  Archive for Mathematical Logic {\bf 55} (2016), nos. 3-4, 562?592.

\bibitem[Vasc]{vaseyc}
Sebastien Vasey, \emph{Building independence relations in abstract elementary classes},  Annals of Pure and Applied Logic {\bf 167} (2016), no. 11, 1029?1092.
\bibitem[Van06]{vandierennomax}
Monica VanDieren, \emph{Categoricity in abstract elementary classes with no
  maximal models}, Annals of Pure and Applied Logic \textbf{141} (2006),
  108--147.
  
  \bibitem[Van13]{nomaxerrata}
Monica VanDieren, \emph{Erratum to ``Categoricity in abstract elementary classes with no
  maximal models''}, Annals of Pure and Applied Logic \textbf{164} (2013),
  131--133.

\bibitem[Vi06]{Vi} 
Andres Villaveces (moderator), \emph{Abstract Elementary Classes workshop open problems session}, Summer 2006, \verb+http://www.aimath.org/WWN/categoricity/categoricity.pdf+


\end{thebibliography}
%\bibliographystyle{amsalpha}

\end{document}